\documentclass{article}

\RequirePackage[OT1]{fontenc}
\RequirePackage{amsthm,amsmath}
\RequirePackage[colorlinks,citecolor=blue,urlcolor=blue]{hyperref}

\usepackage{amsmath}%
\usepackage{amsfonts}%
\usepackage{amssymb}%

\numberwithin{equation}{section}
\theoremstyle{plain}
\newtheorem{theorem}{Theorem}[section]
\newtheorem{lemma}{Lemma}[section]
\newtheorem{corollary}{Corollary}[section]
\theoremstyle{remark}

\def \R  {\mathbb{R}}

\def \P  {\mathbb{P}}
\def \E  {\mathbb{E}}

\def \pgd {LDP }

\def \de {d_{\Sigma}}
\def \dg {d_{e}}

\def \cvloi {\stackrel {w}{\rightarrow}}
\def \cvw1 {\stackrel {W_1}{\rightarrow}}

\def \nin {\not\in}

\def \Cb {C_{b}(\mathbb{R}_+ \times \Sigma)}

\begin{document}

\title{Large deviations for bootstrapped  empirical measures}

\author{ Jos\'e Trashorras\footnote{\texttt{xose@ceremade.dauphine.fr}}\quad and Olivier Wintenberger\footnote{\texttt{wintenberger@ceremade}}\\\small
 Universit\'e Paris-Dauphine, Ceremade, \\\small
 Place du Mar\'echal de Lattre de Tassigny, \\\small
 75775 Paris Cedex 16 France. }

\maketitle




\begin{abstract}
We investigate the Large Deviations properties of bootstrapped empirical 
measure with exchangeable weights. Our main result shows in great generality how the resulting rate function
combines the LD properties of both the sample weights and the observations.
As an application we recover known conditional and unconditional LDPs and obtain some new ones.
\end{abstract}




\section{Introduction and main results }


We say that a sequence of Borel probability measures $(P^{n})_{n \geq 1}$ on a topological space 
$\mathcal{Y}$ obeys a Large Deviation Principle (hereafter abbreviated LDP)
with rate function $I$ if $I$ is a non-negative, lower semi-continuous function defined on $\mathcal{Y}$ 
such that

$$ 
- \inf_{y \in A^{o}} I(y) \leq \liminf_{n \rightarrow \infty} \frac{1}{n} \log P^{n}(A)
     \leq \limsup_{n \rightarrow \infty} \frac{1}{n} \log P^{n}(A) \leq -
     \inf_{y \in \bar{A}} I(y)
$$

\noindent
for any measurable set $A \subset \mathcal{Y}$, whose interior is denoted by $A^{o}$ and closure by $\bar{A}$. 
If the level sets $\{y : I(y) \leq \alpha \}$ are compact for every $\alpha < \infty$, $I$ is called 
a good rate function. With a slight abuse of language we say that a sequence of random variables obeys an 
\pgd when the sequence of measures induced by these random variables obeys an
LDP. For a background on the theory of large  deviations, see Dembo and Zeitouni \cite{dz98} 
and references therein.

\bigskip

Our framework is the following: We are given a triangular array
$((W_{i}^{n})_{1 \leq i \leq n})_{n \geq 1}$ of $\mathbb{R}_+$-valued random variables 
defined on a probability space $(\Omega,\mathcal{A},\mathbb{P})$ and  such that 

\medskip

\begin{itemize}
\item[\bf{(H1)}] For every $n \geq 1$  we have $\sum_{i=1}^n W_{i}^{n} = n$.
\item[\bf{(H2)}] For every $n \geq 1$  the vector $(W_{1}^n,\dots,W_{n}^n)$ is $n$-exchangeable i.e. 
for every element $\sigma$
of  $\mathfrak{S}_n$ the set of permutations of $\{1,\dots,n\}$ the vectors $(W_{1}^n,\dots,W_{n}^n)$
and $(W_{\sigma(1)}^n,\dots,W_{\sigma(n)}^n)$
have the same distribution.
\end{itemize}

\medskip

\noindent
We shall assume that the $W^{n}_i$'s have some LD properties to be detailed below.
We are further given a triangular array $((x^{n}_{i})_{1 \leq i \leq n})_{n \geq 1}$ of elements 
of a Polish space
$(\Sigma, d_{\Sigma})$ such that
\begin{equation} \label{fixedobs}
\mu^{n} =\frac 1n \sum_{i=1}^{n} \delta_{x_{i}^{n}} \cvloi \mu \in M_{1}(\Sigma)
\end{equation}
where $\cvloi$ stands for weak convergence in  the space of Borel probability measures on
$\Sigma$. Let us recall that $\mu^{n} \cvloi \mu$ if and only if
for every real-valued, bounded and continuous application $f$ defined on $\Sigma$
we have $\int_{\Sigma} f(x) \mu^n (dx) \rightarrow \int_{\Sigma} f(x) \mu (dx)$.
Our goal in the present paper 
is to investigate the LD properties of
\begin{equation} \label{object}
\mathcal{L}^n = \frac{1}{n} \sum_{i=1}^n W^{n}_{i} \delta_{x_{i}^{n}}.
\end{equation}
\noindent
So far LD properties of {\it families} of randomly weighted empirical measures 
like (\ref{object}) have been established
in  the particular case of {\it independent and identically distributed} 
$W_{1}^n,\dots,W_{n}^n$, see e.g. \cite{egp93,dz96,n02}.
Results on the LD properties of some particular cases of $\mathcal{L}^n$ as considered here, i.e.
{\it with exchangeable} $W_{1}^n,\dots,W_{n}^n$
are available but their proofs rely on
the definition of the chosen sampling weights, see e.g. \cite{bj88,dz96,ck96}. Hence, the present paper gives the first 
derivation of the LD properties of the family of empirical measures $\mathcal{L}^n$ under the natural 
assumptions $\textbf{(H1-H2)}$. The interest in considering families rather
than particular cases lies on the fact that only an upper level of abstraction can reveal the 
mechanisms that really enter into play. 

\bigskip

Here we shall consider $\mathcal{L}^n$ from the
bootstrap point of view and adopt its vocabulary \cite{Efron1982,h92,bb95,g97}. 
In 1979, in a landmark paper, Efron proposed the following idea:
When given a realization $x^n_1,\dots,x^ n_n$ of random variables $X^n_1,\dots,X^ n_n$ one can easily
obtain "additional data" by sampling independent and $\frac 1n \sum_{i=1}^{n} \delta_{x_{i}^{n}}$-distributed
random variables $X^{\ast}_1,\dots,X^{\ast}_m,\dots$. It amounts to sample with replacement from an urn
which composition is described by $\frac 1n \sum_{i=1}^{n} \delta_{x_{i}^{n}}$. Often, this is 
computationally cheap and theoretical studies are available to assess
the quality of the distribution of $H(\frac{1}{m} \sum_{i=1}^m  \delta_{X_{i}^{\ast}})$ in approximating the distribution of
$H(\frac{1}{n} \sum_{i=1}^n  \delta_{X_{i}^n})$ which make it all worthwhile. A rich literature started flourishing 
on the ground of this idea. It was soon noticed that the preceding procedure is not the only possible one
and that it can be generalized so that  $\frac{1}{n} \sum_{i=1}^n W^{n}_{i} \delta_{x_{i}^{n}}$ under conditions $\textbf{(H1-H2)}$
is the right object to be considered. For example, Efron's bootstrap corresponds to  
$(W_{1}^n,\dots,W_{n}^n)$ distributed according to a Multinomial law. The literature on the subject developed into two
complementary directions:
"conditional" results where $x^n_1,\dots,x^ n_n$ are fixed observations filling some conditions
and $H(\frac{1}{n} \sum_{i=1}^n W^{n}_{i} \delta_{x_{i}^{n}})$ is considered and "unconditional" results where the
$x^n_1,\dots,x^ n_n$ are allowed to fluctuate and $H(\frac{1}{n} \sum_{i=1}^n W^{n}_{i} \delta_{X_{i}^{n}})$
is considered instead.  

Classically the bootstrap scheme is said to be efficient when it mimics  the behavior of 
$H(\frac{1}{n} \sum_{i=1}^n \delta_{X_{i}^{n}})$ and one distinguishes between "conditional efficiency"
and "unconditional efficiency". For example,  Praestgaard and Wellner investigated the Central Limit behavior 
of  $\mathcal{L}^n$ in \cite{pw93}; they gave a necessary and sufficient condition on the second moments properties of the $W_i^n$
which works for both conditional and unconditional efficiency. 
Later, Hall and Mammen studied the efficiency of the bootstrap schemes in the Edgeworth expansions at the second 
order \cite{h94}. There, conditions on the fourth order cumulants of the weights are required. 
In a similar context, Wood already showed in \cite{w00} that Efron's bootstrap is  efficient 
to mimic the empirical mean in the moderate deviations regime
for observations satisfying the Cramer condition but not for heavier tails.

Following Barbe and Bertail \cite{bb95}, we say  that a bootstrap scheme is LD-efficient when the bootstrapped 
empirical measure has the same LD properties as the original empirical measure. It is a very strong property: thinking of  percentile
bootstrap's confidence intervals, LD-efficiency says that the relative coverage accuracy tends to 1 exponentially fast.
As an application of our general approach, we will be able to discuss both conditional and unconditional LD-efficiency
for many classical choices of $((W_{i}^{n})_{1 \leq i \leq n})_{n \geq 1}$ and/or 
$((X_{i}^{n})_{1 \leq i \leq n})_{n \geq 1}$. Actually we go further in the sense that as an application of our approach we 
obtain LD results that are new in the literature like
e.g. an unconditional LDP for Efron's bootstrap and conditional and unconditional LDP's for iid weighted bootstrap and $k$-blocks bootstraps.

\bigskip

Now let us  describe our results more precisely. To this end we need to introduce some more notations.
We consider 

$$M_{1}^{1} (\mathbb{R}_+) = \left\{ \rho \in M_{1}(\mathbb{R}_+) \ : \ \int_{\mathbb{R}_+} x \ \rho(dx) = 1 \right\}$$
a subset of
$$\mathcal{W}^{1} (\mathbb{R}_+) = \left\{ \rho \in M_{1}(\mathbb{R}_+) \ : \ \int_{\mathbb{R}_+} x \ 
\rho(dx) < \infty \right\}$$
the so-called Wasserstein space of order 1 on $\mathbb{R}_+$.
Both $M_{1}^{1} (\mathbb{R}_+)$ and $\mathcal{W}^{1} (\mathbb{R}_+)$ are endowed 
with the Wasserstein-1 distance

$$W_{1,|\cdot|}(\rho,\gamma ) = \inf_{\pi \in \mathcal{C}(\rho,\gamma )} 
\left\{ \int_{\mathbb{R}_+ \times \mathbb{R}_+}
|x-y| \pi(dx,dy) \right\}$$
where $\mathcal{C}(\rho,\gamma )$ is the subset of $M_{1}(\mathbb{R}_+ \times \mathbb{R}_+)$ of 
couplings of $\rho$ and $\gamma $
i.e. the set of Borel probability measures $\pi$ on $\mathbb{R}_+ \times \mathbb{R}_+$ such that their first marginal
$\pi_1$ is $\rho$ and second marginal $\pi_2$ is $\gamma $. 
In addition to \textbf{(H1-H2)} we shall assume that the $((W_{i}^{n})_{1 \leq i \leq n})_{n \geq 1}$ are 
such that

\medskip

\begin{itemize}
\item[\bf{(H3)}] The sequence $(\mathcal{S}^n = \frac{1}{n} \sum_{i=1}^n \delta_{W^{n}_i})_{n \geq 1}$ 
satisfies an LDP on
$M_{1}^{1} (\mathbb{R}_+)$ endowed with the $W_{1,|\cdot|}$ distance with good rate function $I^W$.
\end{itemize} 

\medskip

\noindent
In Section \ref{catalogue} we prove that \textbf{(H3)} is valid for a broad collection of sampling weights.
Actually we shall see that the LD properties of $(\mathcal{L}^n)_{n \geq 1}$ directly follow 
from the LD properties of 
$$\mathcal{V}^n = \frac{1}{n} \sum_{i=1}^n \delta_{(W^{n}_{i}, x_{i}^{n})}$$
on the set
$$\mathcal{M}_{1}^{1} (\mathbb{R}_+ \times \Sigma) = \left\{ \rho \in M_{1} (\mathbb{R}_+ \times \Sigma) \ : \ 
\int_{\mathbb{R}_+} x \rho_1(dx)  = 1  \right\}$$
endowed with the distance 
$$\Delta(\rho,\gamma ) = W_{1,|\cdot|}(\rho_1,\gamma_1) + \beta_{BL,\delta}(\rho,\gamma )$$
where
$$\beta_{BL,\delta}(\rho,\gamma) = \sup_{f \in C_b (\mathbb{R}_+ \times \Sigma) \atop
||f||_{\infty} + ||f||_{L,\delta} \leq 1} \left\{ \left| \int_{\mathbb{R}_+} f(x) \rho(dx) -
\int_{\mathbb{R}_+} f(x) \gamma (dx) \right| \right\}$$
is the so-called dual-bounded Lipschitz metric on $M_{1} (\mathbb{R}_+ \times \Sigma)$ and as usual

$$||f||_{\infty} = \sup_{x \in \mathbb{R}_+ \times \Sigma} |f(x)|, \ \ \ \ \ \ 
||f||_{L,\delta} = \sup_{x,y \in \mathbb{R}_+ \times \Sigma \atop x \neq y} \frac{|f(x) - f(y)|}{\delta(x,y)},$$
$\delta$ is a metric on $\mathbb{R}_+ \times \Sigma$ compatible with the product topology
and $C_b (\mathbb{R}_+ \times \Sigma)$ is the set of real-valued, bounded and continuous applications defined on
$\mathbb{R}_+ \times \Sigma$.
For any two probabilities $\rho,\nu$ on a measurable space $(E,\mathcal{E})$ we denote by

$$
H(\nu| \rho) = \left\{
\begin{array}{cl}
\int_E {\mbox d}\nu \log \frac{{\mbox d}\nu}{{\mbox d}\rho} & \mbox{if \ } \nu \ll \rho \\
+ \infty & \mbox{otherwise}
\end{array}
\right.
$$
the relative entropy of $\nu$ with respect to $\rho$. To any $\rho(dw,dx) \in M_1(\mathbb{R}_+ \times \Sigma)$
we associate $\rho_x(dw) \in M_1(\R_+)$ (resp. $\rho_w(dx) \in M_1(\Sigma)$) a stochastic kernel 
which is the conditional distribution of the first (resp. second) marginal of $\rho$ given the second (resp. first).
We summarize this by $\rho(dw,dx) = \rho_x (dw) \otimes \rho_2(dx)$ (resp. $\rho(dw,dx) = \rho_1(dw) \otimes \rho_w (dx)$).
If $\nu, \gamma \in M_1(\mathbb{R}_+ \times \Sigma)$ are such that $\nu_1=\gamma_1=\theta$ (resp. $\nu_2=\gamma_2=\theta$)
then 
\begin{equation}
\label{tren}
H(\nu|\gamma) = \int_{\R_+} H(\nu_w |\gamma_w ) \theta(dw)
\end{equation}
(resp. $H(\nu|\gamma) = \int_{\Sigma} H(\nu_x |\gamma_x ) \theta(dx)$), see Lemma 1.4.3 in \cite{de97}. Our first result is the following

\begin{theorem} \label{th1}
The sequence 
$(\mathcal{V}^n)_{n \geq 1}$ satisfies an LDP on $\mathcal{M}^{1}_{1}(\mathbb{R}_+ \times \Sigma)$
endowed with the distance $\Delta$ with good rate function

$$
\mathcal{J}(\rho ; \mu) =  \left\{
\begin{array}{cl}
H(\rho| \rho_1 \otimes \mu) + I^{W}(\rho_1) & \mbox{if \ } \rho_2 = \mu  \\
+ \infty & \mbox{otherwise}.
\end{array}
\right.
$$
\end{theorem}
\noindent
The reason for working on $\mathcal{M}^{1}_{1}(\mathbb{R}_+ \times \Sigma)$ endowed with the distance $\Delta$
is to make the following map 
$$
\begin{array}{rccc}
F: & \mathcal{M}^{1}_{1}(\mathbb{R}_+ \times \Sigma) & \rightarrow & M_1(\Sigma) \\
   &  \rho(dw,dx) & \mapsto & \int_{\mathbb{R}_+} w \rho(dw,dx)
\end{array}
$$
well-defined and continuous once $M_1(\Sigma)$ is endowed with the weak convergence topology.
By contraction (see Theorem 4.2.1 in \cite{dz98}) an LDP for $(\mathcal{L}^n)_{n \geq 1}$ easily follows
from Theorem \ref{th1}
\begin{corollary} \label{th2}
The sequence 
$(\mathcal{L}^{n})_{n \geq 1}$ satisfies an LDP on $M_{1}(\Sigma)$
endowed with the weak convergence topology with good rate function
\begin{eqnarray*}
\mathcal{K}(\nu ; \mu) & =  & \inf_{\rho : F(\rho)=\nu} \mathcal{J}(\rho ; \mu) \\
& = &  \inf_{\rho_x : F(\rho_x \otimes \mu)=\nu} \left\{ \int_{\Sigma} H(\rho_x | \rho_1) \mu(dx) + I^W(\rho_1) \right\}.
\end{eqnarray*}
\end{corollary}
\noindent
We shall often see in applications that $I^W(\nu)=H(\nu|\xi)$ for some $\xi \in M_1(\R_+)$ such that
\begin{equation} \label{douane}
 \Lambda_{\xi}(\alpha) = \log \int_{\mathbb{R}_+} e^{\alpha x} \xi(dx) < \infty
\end{equation}
\noindent
for every $\alpha \in \R$. We denote by 
\begin{equation} \label{francoisimenhaus}
\Lambda_{\xi}^{\ast}(x) = \sup_{\alpha \in \R} \left\{ \alpha x -  \Lambda_{\xi}(\alpha) \right\}. 
\end{equation}
\noindent
the Fenchel-Legendre transform of $\Lambda_{\xi}$. In this particular case we obtain a generalization of Kullback inequality 
\begin{lemma} \label{clement}
If $I^W(\nu)=H(\nu|\xi)$ for some $\xi \in M_1(\R_+)$ then for every $\nu,\mu \in M_1(\Sigma)$ we have
\begin{equation} \label{varlin}
\mathcal{K}(\nu ; \mu) \geq \left\{
\begin{array}{cl}   
\int_{\Sigma}  \Lambda_{\xi}^{\ast}( \frac{d\nu}{d\mu}(x) ) \mu(dx) & \mbox{ if } \nu \ll \mu \\
+ \infty & \mbox{ otherwise. } \\
\end{array}
\right.
\end{equation}
Moreover, assuming that (\ref{douane}) holds
for every $\alpha \in \R$, the preceding inequality turns out to be an equality for every $\nu,\mu \in M_1(\Sigma)$ if and only if 
\begin{equation} \label{gargamel}
\lim_{\alpha \rightarrow - \infty} \Lambda_{\xi}^{'}(\alpha) = 0 \mbox{   and   }
\lim_{\alpha \rightarrow + \infty} \Lambda_{\xi}^{'}(\alpha) = + \infty.
\end{equation}
\end{lemma}
\noindent
Remark that the right hand side of the inequality \eqref{varlin} is the rate function obtained in case of iid weights in \cite{dz96}. 

\noindent

Next we allow the $x^{n}_{i}$'s to fluctuate and consider a triangular array 
$((X_{i}^{n})_{1 \leq i \leq n})_{n \geq 1}$ of $\Sigma$-valued random variables 
defined on $(\Omega,\mathcal{A},\mathbb{P})$ such that

\medskip

\begin{itemize}
\item[\bf{(H4)}] The sequence $(\mathcal{O}^n = \frac{1}{n} \sum_{i=1}^n \delta_{X^{n}_i})_{n \geq 1}$ satisfies 
an LDP on
$M_{1} (\Sigma)$ endowed with the weak convergence topology with good rate function $I^X$.
\item[\bf{(H5)}] For every $n \geq 1$ the vectors $(X_{1}^{n},\dots,X_{n}^{n})$ 
and $(W_{1}^{n},\dots,W_{n}^{n})$ are independent.
\end{itemize}

\medskip

\noindent
An LDP for 
$$V^n = \frac{1}{n} \sum_{i=1}^n \delta_{(W^{n}_{i}, X_{i}^{n})}$$
holds as a consequence of Theorem \ref{th1} and Theorem 2.3 in \cite{g96}

\begin{theorem} \label{th3}
The sequence 
$(V^n)_{n \geq 1}$ satisfies an LDP on $\mathcal{M}^{1}_{1}(\mathbb{R}_+ \times \Sigma)$
endowed with the distance $\Delta$ with good rate function

\begin{equation} \label{tarantino}
J(\rho) =  H(\rho| \rho_1 \otimes \rho_2) + I^{W}(\rho_1) + I^{X}(\rho_2).
\end{equation}
\end{theorem}
\noindent
Again, by contraction, an LDP for 
$$L^n = \frac{1}{n} \sum_{i=1}^n W^{n}_i \delta_{X_{i}^{n}}$$
easily follows from Theorem \ref{th3}
\begin{corollary} \label{th4}
The sequence 
$(L^{n})_{n \geq 1}$ satisfies an LDP on $M_{1}(\Sigma)$
endowed with the weak convergence topology with good rate function

\begin{eqnarray}
K(\nu) & =  & \inf_{\rho : F(\rho)=\nu} J(\rho) \nonumber \\
& = & \inf_{\rho_2 \in M_1(\Sigma)} \left\{ \inf_{\rho_x : F(\rho_x \otimes \rho_2)=\nu} 
\left\{ \int_{\Sigma} H(\rho_x | \rho_1) \rho_2(dx) + I^W(\rho_1) \right\} + I^X(\rho_2) \right\} \nonumber \\
& = & \inf_{\rho_2 \in M_1(\Sigma)} \left\{ \mathcal{K}(\nu ; \rho_2) + I^X(\rho_2) \right\}. \label{shrek}
\end{eqnarray}
It follows that for every $\nu \in M_{1}(\Sigma)$ we have $K(\nu) \leq I^{X}(\nu)$.
\end{corollary}

\noindent
The latter inequality illustrates the \textit{smoothing effect} of the random weights $W^n_1,\dots,W^n_n$
on the distribution of $L^n$. We shall see in most examples that for classical choices of $((W_{i}^{n})_{1 \leq i \leq n})_{n \geq 1}$ 
and/or $((X_{i}^{n})_{1 \leq i \leq n})_{n \geq 1}$ there exists at least one $\nu \in M_1(\Sigma)$ such that
$K(\nu) < I^X(\nu)$. Nevertheless, we have the following

\begin{corollary} \label{contador} A necessary and sufficient condition on $((W_{i}^{n})_{1 \leq i \leq n})_{n \geq 1}$ to ensure
that for every $((X_{i}^{n})_{1 \leq i \leq n})_{n \geq 1}$ and every $\nu \in M_1(\Sigma)$ we have
$K(\nu) = I^X(\nu)$ is that for every $\nu,\zeta \in M_1(\Sigma)$ 
$$
\mathcal{K}(\nu;\zeta) = \left\{
\begin{array}{cl}
0 & \mbox{ if } \nu= \zeta \\
+ \infty & \mbox{ otherwise. }
\end{array}
\right.
$$
\end{corollary}

\noindent
It follows from Lemma \ref{clement} that is condition is satisfied e.g. when
the sampling weights are such that $I^W(\nu)=H(\nu|\delta_1)$ as for 
the delete-$h(n)$ jacknife with $h(n)=o(n)$, see Corollary \ref{cavendish} below.
 
\bigskip

The paper is structured as follows. Section 2 presents several applications of our main results. There we discuss
efficiency issues for both conditional and unconditional LDPs. The rest of the paper is devoted to proofs.
In particular Section \ref{catalogue} contains the proofs of several sampling weights LDPs. Some of these are new 
to the literature.


\section{Examples of applications} \label{anigo}


\noindent
In this section we investigate the LD properties of $(\mathcal{L}^n)_{n \geq 1}$ and $(L^n)_{n \geq 1}$
for several particular choices of $((W_{i}^{n})_{1 \leq i \leq n})_{n \geq 1}$ and/or 
$((X_{i}^{n})_{1 \leq i \leq n})_{n \geq 1}$. To specify our results we need some more notations.
For every $\lambda > 0$ we shall denote by $\mathcal{P}(\lambda)$ the Poisson distribution with parameter $\lambda$
and by $\mathcal{Q}(\lambda)$ the distribution of a random variable $Y$ such that $\lambda Y$ 
is $\mathcal{P}(\lambda)$-distributed. More generally, for every $\lambda, \gamma > 0$ we shall denote 
by  $\mathcal{F}(\lambda,\gamma)$ the distribution of a random variable $Y$ such that $\lambda Y$ 
is $\mathcal{P}(\gamma)$-distributed.
For every positive integers $m$ and $n$ and every $n$-tuple of non-negative numbers $(p^n_1,\dots,p^n_n)$
such that $\sum_{i=1}^n p^n_i =1$ we shall denote by Mult$_{n}(m,(p^n_1,\dots,p^n_n))$ the distribution
of $(Y_1,\dots,Y_n)$ the numbers of balls found in $n$ urns labeled $1,\dots,n$ when $m$ balls are thrown in these urns 
independently, each having probability $p^n_1$ to fall in the urn labeled $1$, 
probability $p^n_2$ to fall in the urn labeled $2$, etc...

\subsection{Efron's bootstrap and "$m$ out of $n$" bootstrap}
\label{fronfron}
\noindent
For every $m,n \geq 1$ the weights $(W^n_1,\dots,W^n_n)$ for the "$m$ out of $n$" bootstrap
are defined such that  
$$
\frac{m}{n}(W^n_1,\dots,W^n_n)\sim \mbox{Mult}_{n}(m,(1/n,\dots,1/n)).
$$ 
Classical Efron's bootstrap corresponds to $m=n$. It emerges from sampling
with replacement from the urn containing the observed data.
We shall assume that $m=m(n)$ and that the sequence $(\lambda_n = m(n)/n)_{n \geq 1}$ satisfies 
$\lim_{n \rightarrow \infty} \lambda_n = \lambda > 0$. Quite surprisingly we could not find in the literature 
a reference for the following, even in the simple $m(n)=n$ case

\begin{theorem} \label{lemmaE} The sequence $(\mathcal{S}^n = \frac{1}{n} \sum_{i=1}^n  \delta_{W^n_i})_{n \geq 1}$ 
obeys an LDP
on $M^1_1(\R_+)$ endowed with the $W_{1,|\cdot|}$ distance  with good rate function $H(\cdot|\mathcal{Q}(\lambda))$.
\end{theorem}

\noindent
First we consider fixed observations $((x^n_i)_{1 \leq i \leq n})_{n \geq 1}$ such that (\ref{fixedobs}) holds.

\begin{corollary} \label{nkoulou} The sequence 
$(\mathcal{L}^{n})_{n \geq 1}$ satisfies an LDP on $M_{1}(\Sigma)$
endowed with the weak convergence topology with good rate function
$$
\mathcal{K}(\nu;\mu) = \lambda H(\nu|\mu).
$$
\end{corollary}

\noindent
By properly rescaling we immediately obtain for every $\lambda >0$ and every measurable $A \subset M_1(\Sigma)$ that
 
\begin{multline}
- \inf_{\nu \in A^{o}} H(\nu|\mu) \leq \liminf_{n \rightarrow \infty} \frac{1}{m(n)} \log \P(\mathcal{L}^{n} \in A)  \leq \\
 \leq \limsup_{n \rightarrow \infty} \frac{1}{m(n)} \log \P(\mathcal{L}^{n} \in A) \leq -
     \inf_{\nu \in \bar{A}} H(\nu|\mu). \label{scaling}
\end{multline}

\noindent
Next we investigate the LD 
properties of $(L^n)_{n \geq 1}$ without any other assumption 
that $\textbf{(H4-H5)}$. To this end we introduce

$$\mathcal{Z} = \left\{ \eta \in M_1(\Sigma) \ : \ I^X(\eta) = 0 \right\}.$$
It follows from Corollary \ref{th4} and \ref{nkoulou} that

\begin{corollary} \label{amalfitano} The sequence 
$(L^{n})_{n \geq 1}$ satisfies an LDP on $M_{1}(\Sigma)$
endowed with the weak convergence topology with good rate function $K$ such that
$$
K(\nu) = \inf_{\zeta \in M_1(\Sigma)} \{ \lambda H(\nu| \zeta) + I^X(\zeta) \} 
\leq \lambda \inf_{\eta \in \mathcal{Z}}  H(\nu|\eta).
$$
\end{corollary}

\noindent
Now we consider some particular cases for 
$((X^n_i)_{1 \leq i \leq n})_{n \geq 1}$.
First, we assume that 
for every $n \geq 1$ the random variables $X^n_1,\dots,X^n_n$ are independent and identically $\mu$-distributed. Then
$\frac{1}{n} \sum_{i=1}^n \delta_{X^n_i} \cvloi \mu \ a.s.$ (see Theorem 11.4.1 in \cite{d02}) 
and Corollary \ref{nkoulou} can be interpreted
 as a conditional LDP.
Hence, any "$m$ out of $n$" bootstrap such that 
$\lim_{n \rightarrow \infty} m(n)/n = 1$ (in particular Efron's bootstrap) leads to a conditional LDP 
that coincides with the original LDP in this case. This was first established in \cite{bj88} for Efron's Bootstrap
and in \cite{ck96} in the general case. Actually the $X^n_1,\dots,X^n_n$ need not be iid
but only that the associated empirical measures satisfy an LDP with rate function $H(\cdot|\mu)$ for Efron's bootstrap to be
conditionally LD-efficient.
Corollary \ref{amalfitano} completes the previous result with an unconditional LDP. In this particular case
$I^X(\zeta)=H(\zeta|\mu)$ and by taking e.g. $\mu=\frac{9}{10}\delta_0 + \frac{1}{10}\delta_1$,
$\nu=\frac{1}{10}\delta_0 + \frac{9}{10}\delta_1$, $\zeta=\frac{1}{2}\delta_0 + \frac{1}{2}\delta_1$ and
$\lambda=1$ we observe that $K(\nu)\leq H(\nu|\zeta) + H(\zeta|\mu) < H(\nu|\mu)$ hence Efron's bootstrap is not unconditionally LD-efficient.
Straightforward use of the same kind of arguments shows that this remark also holds true when $X^n_1,\dots,X^n_n$
is the result of sampling without replacement from an urn with suitable properties (see Theorem 7.2 in \cite{dz98} for the reference LDP)
or when $X^n_1,\dots,X^n_n$ are the $n$ first components of an infinitely exchangeable sequence of random variables 
(see \cite{dz92} for the reference LDP).

\subsection{Iid weighted bootstrap}

The weights $(W_1^n,\dots,W^n_n)$ for an iid-weighted bootstrap are defined on the ground 
of a sequence $Y_1,\dots,Y_n,\dots$
of $\mathbb{R}_+$-valued independent random variables with common distribution $\xi$. We shall assume that 
for every $\alpha > 0$ we have 
$\Lambda_{\xi}(\alpha)  < \infty$ 
and that $\Lambda_{\xi}^{\ast}(0)=\infty$ (or equivalently $\P(Y_1=0)=0$).
The weights $(W_1^n,\dots,W^n_n)$ are defined by

$$W^n_1 = \frac{Y_1}{\frac{1}{n} \sum_{i=1}^n Y_i},\dots,W^n_i = \frac{Y_i}{\frac{1}{n} \sum_{i=1}^n Y_i},
 \dots,W^n_n = \frac{Y_n}{\frac{1}{n} \sum_{i=1}^n Y_i}.$$ 
In order to describe the LD behavior of  $(\mathcal{S}^n = \frac{1}{n} \sum_{i=1}^n \delta_{W^{n}_{i}})_{n \geq 1}$ 
we introduce the  map
$$
\begin{array}{rccc}
\mathcal{G}: & \mathcal{W}^{1}(\mathbb{R}_+) \times \mathbb{R}_{+}^{*} & \rightarrow & \mathcal{W}^{1}(\mathbb{R}_+) \\
   &  (\rho,m) & \mapsto & \mathcal{G}(\rho,m): A \in \mathcal{B}_{\mathbb{R}_+}  \mapsto \rho(m A)
\end{array}
$$
where for every Borel set $A \in \mathcal{B}_{\mathbb{R}_+}$ and every $m>0$ we write
$$mA = \left\{ x \in \mathbb{R}_+, \exists y \in A : x=my \right\}.$$

\noindent
The continuity of $\mathcal{G}$ is the main argument in the proof of the following
\begin{theorem} \label{pgdiid}
The sequence $(\mathcal{S}^n = \frac{1}{n} \sum_{i=1}^n \delta_{W^{n}_{i}})_{n \geq 1}$ 
satisfies an LDP on
$M^1_1 (\mathbb{R}_+)$ endowed with the $W_{1,|\cdot|}$ distance with good rate function 
$$
I^W(\rho) = \inf_{m > 0} \left\{ H( \mathcal{G}(\frac{1}{m},\rho )| \xi ) \right\}.
$$
\end{theorem}
\noindent
The previous LDP leads to

\begin{corollary} \label{vockler} For every $\nu,\mu \in M_1(\Sigma)$ we have
$$
\mathcal{K}(\nu ; \mu) \geq \left\{
\begin{array}{cl}   
\inf_{m > 0} \int_{\Sigma}  \Lambda_{\xi}^{\ast}( m \frac{d\nu}{d\mu}(x) ) \mu(dx) & \mbox{ if } \nu \ll \mu \\
+ \infty & \mbox{ otherwise. } \\
\end{array}
\right.
$$
Moreover, the preceding turns out to be an equality for every $\nu,\mu \in M_1(\Sigma)$ if and only if
$$ 
\lim_{\alpha \rightarrow - \infty} \Lambda_{\xi}^{'}(\alpha) = 0 \mbox{   and   }
\lim_{\alpha \rightarrow + \infty} \Lambda_{\xi}^{'}(\alpha) = + \infty.
$$
\end{corollary}

\noindent
It follows from the previous corollary that there is no distribution $\xi$ such that for every
$\nu,\mu \in M_1(\Sigma)$ the identity $\mathcal{K}(\nu;\mu) =  H(\nu | \mu)$ holds.
Indeed, as soon as there exists $\nu,\mu \in M_1(\Sigma)$ such that
$\frac{d\nu}{d\mu}(x) = 0$ on a set $A$ such that $\mu(A) > 0$ one has $\mathcal{K}(\nu;\mu)=\infty$
while it could be possible that $H(\nu | \mu) < \infty$. In words there is no choice of $\xi$ for which one
 gets a conditional LDP that coincides with the original one for 
$X^n_1,\dots,X^n_n$ independent and $\mu$-distributed. It is clearly due to the fact that 
$\Lambda_{\xi}^{\ast}(0)=\infty$ forces all the weights $W^n_1,\dots,W^n_n$ to be positive which is to be compared to e.g.
Efron's bootstrap. Finally, as for Efron's bootstrap, one can construct examples to show that in most classical cases
the iid-bootstrap is not unconditionally LD-efficient. 

\subsection{The multivariate hypergeometric bootstrap}

Let $K$ be a fixed integer number such that $K \geq 2$. 
The multivariate hypergeometric bootstrap emerges from the following urn scheme: Put $K$ copies
of each observed data in an urn so that the urn contains $Kn$ elements then draw from this urn a sample of size $n$
without replacement. The sampling weights $(W^{n}_{1},\dots,W^{n}_{n})$ take their values in $\{0,1,\dots,K\}$ under the constraint
$\sum_{i=1}^n W^{n}_{i} =n$ and are distributed according to

$$\mathbb{P}(W^{n}_{1}=w^{n}_{1},\dots,W^{n}_{n}=w^{n}_{n}) = \frac{C_K^{w^{n}_{1}} \cdots C_K^{w^{n}_{n}}}{C_{nK}^n}.$$
Let us denote by $\mathfrak{B}(K,\frac{1}{K})$ the Binomial distribution with parameters $K$ and $\frac{1}{K}$.

\begin{theorem} \label{binopgd} The sequence $(\mathcal{S}^n = \frac{1}{n} \sum_{i=1}^n \delta_{W^{n}_{i}})_{n \geq 1}$ satisfies 
an LDP on
$M^1_1 (\mathbb{R}_+)$ endowed with the $W_{1,|\cdot|}$ distance with good rate function 

$$
I^W(\rho) = 
H(\rho| \mathfrak{B}(K,\frac{1}{K})) 
$$
\end{theorem}
\noindent
Consider fixed observations $((x^n_i)_{1 \leq i \leq n})_{n \geq 1}$ such that (\ref{fixedobs}) holds, we obtain
\begin{corollary}The sequence 
$(\mathcal{L}^{n})_{n \geq 1}$ satisfies an LDP on $M_{1}(\Sigma)$
endowed with the weak convergence topology with good rate function
$$
\mathcal{K}(\nu ; \mu) \geq \left\{
\begin{array}{cl}   
 \int_{\Sigma}  \Lambda_{\mathfrak{B}(K,\frac{1}{K})}^{\ast}( \frac{d\nu}{d\mu}(x) ) \mu(dx) & \mbox{ if } \nu \ll \mu \\
+ \infty & \mbox{ otherwise. } \\
\end{array}
\right.
$$
\end{corollary}
\noindent
We only obtain an inequality since $\mathfrak{B}(K,\frac{1}{K})$ does not satisfy condition (\ref{gargamel}).
Again, there is no integer $K$ such that for every
$\nu,\mu \in M_1(\Sigma)$ the identity $\mathcal{K}(\nu;\mu) =  H(\nu | \mu)$ holds.
Indeed, as soon as there exists $\nu,\mu \in M_1(\Sigma)$ such that
$\frac{d\nu}{d\mu}(x) > K $ on a set $A$ such that $\mu(A) > 0$ one has $\mathcal{K}(\nu;\mu)=\infty$
while it could be possible that $H(\nu | \mu) < \infty$.
Thus all multivariate hypergeometric bootstraps fail to be conditionally LD-efficients for iid observations.
One can construct examples to show that in most classical cases
the multivariate hypergeometric bootstrap fails to be unconditionally LD-efficient.

\subsection{A bootstrap generated from deterministic weights}

The weights for bootstrap schemes defined from deterministic weights are given by 
$$(W^n_1,\dots,W^n_n) =  (w^n_{\sigma_n(1)},\dots,w^n_{\sigma_n(n)})$$ 
where for every $n \geq 1$ 
the $w^n_1,\dots,w^n_n$ are fixed non-negative real numbers such that $\sum_{i=1}^n w^n_i = n$ and
$$\frac{1}{n} \sum_{i=1}^n \delta_{w^n_i} \cvw1 \gamma \in M_1^1(\mathbb{R}_+)$$ 
and $\sigma_n$ is an uniformly over $\mathfrak{S}_n $ distributed random variable. We clearly have
\begin{theorem} The sequence $(\mathcal{S}^n = \frac{1}{n} \sum_{i=1}^n \delta_{W^{n}_{i}})_{n \geq 1}$ satisfies 
an LDP on
$M^1_1 (\mathbb{R}_+)$ endowed with the $W_{1,|\cdot|}$ distance with good rate function 

$$
I^W(\rho) = \left\{
\begin{array}{cl}
0 & \mbox{if \ } \rho = \gamma \\
+ \infty & \mbox{otherwise}.
\end{array}
\right.
$$
\end{theorem}

\noindent
An important special case is the grouped, or delete-$h$ jacknife. The grouped jacknife with group block size $h$
may be viewed as a bootstrap generated by permuting the deterministic weights
$$
\begin{array}{cccc}
(w^n_{1},\dots,w^n_{n}) & = & \left(\underbrace{ \frac{n}{n-h},\dots,\frac{n}{n-h}} , \right. & \left. \underbrace{0,\dots,0} \right) \\
 & & n-h & h
 \end{array}
$$ 
We shall take $h=h(n)$ such that $\lim_{n \rightarrow \infty} h(n)/n = \alpha \in [0,1)$ so
$$\gamma = (1-\alpha) \delta_{\frac{1}{1-\alpha}} + \alpha \delta_0.$$

\begin{corollary} \label{cavendish} If $\alpha > 0$ the sequence 
$(\mathcal{L}^{n})_{n \geq 1}$ satisfies an LDP on $M_{1}(\Sigma)$
endowed with the weak convergence topology with good rate function
$$
\mathcal{K}(\nu;\mu) = \left\{
\begin{array}{cl}
(1-\alpha)H(\nu|\mu) + \alpha H \left( \frac{\mu - (1 -  \alpha) \nu}{\alpha} |\mu \right ) 
& \mbox{  if   } \frac{\mu - (1 -  \alpha) \nu}{\alpha} \in M_1(\Sigma) \\
+ \infty & \mbox{otherwise.}
\end{array}
\right.
$$
If $\alpha =0 $  the sequence 
$(\mathcal{L}^{n})_{n \geq 1}$ satisfies an LDP on $M_{1}(\Sigma)$
endowed with the weak convergence topology with good rate function 
$$
\mathcal{K}(\nu;\mu) = \left\{
\begin{array}{cl}
0 
& \mbox{  if   } \nu=\mu \\
+ \infty & \mbox{otherwise.}
\end{array}
\right.
$$
\end{corollary}

\noindent
Naturally this result coincides with Theorem 7.2.1 in \cite{dz98}. Combining Corollary \ref{th4} and Corollary \ref{cavendish} we 
obtain an unconditional version of the latter result.
To every $\nu \in M_1(\Sigma)$ we associate

$$
\mathcal{E}_{\nu} = \{ \zeta \in M_1(\Sigma) \ : \ \frac{\zeta - (1 -  \alpha) \nu}{\alpha} \in M_1(\Sigma) \}
$$
\begin{corollary} \label{hinault} If $\alpha > 0$ the sequence 
$(L^{n})_{n \geq 1}$ satisfies an LDP on $M_{1}(\Sigma)$
endowed with the weak convergence topology with good rate function $ K(\nu) = \inf_{\zeta \in  
\mathcal{E}_{\nu}} \mathcal{U}(\nu,\zeta)$
where
$$
\mathcal{U}(\nu,\zeta) = 
(1-\alpha)H(\nu|\zeta) + \alpha H \left( \frac{\zeta - (1 -  \alpha) \nu}{\alpha} |\zeta \right ) + I^X(\zeta).
$$
If $\alpha =0 $  the sequence 
$(L^{n})_{n \geq 1}$ satisfies an LDP on $M_{1}(\Sigma)$
endowed with the weak convergence topology with good rate function $ {K}(\nu )=I^X(\nu)$.
\end{corollary}

\subsection{The $k$-blocks bootstraps}

Let us consider the (moving or circular) $k$-block bootstrap. Consider weights from the "$m=n/k$ out of n" bootstrap:
$$
\frac{1}{k} (W^n_1,\dots,W^n_n)\sim \mbox{Mult}_{n}(m,(1/n,\dots,1/n)).
$$ 
The $k$-blocks bootstrapped empirical measure can be written as in \cite{h04} with the formula
$$
\widetilde{\mathcal L}_n=\frac1n\sum_{i=1}^nW_i^n \frac1k \sum_{j\sim i}\delta_{x^n_j}
$$
where $j\sim i$ means that the $ j$ belong to block $i$. For the moving $k$-blocks bootstrap, $j\sim i$ if 
$j \in \{i - k/2,\cdots,i+k/2\}$ miodulo $n$. With slight modifications we could also consider 
the circular $k$-blocks bootstrap where $j\sim i$ if $j\in \{i,\ldots,i+k-1\}$ modulo $n$. Both schemes are asymptotically equivalent 
as soon as $k$ is fixed as it is the case here. Notice that
$$
\widetilde{\mathcal L}_n=\frac1n\sum_{i=1}^n \widetilde{W_i^n} \delta_{x^n_i}
\ \ \ \mbox{with} \ \ \ 
\widetilde{W_i^n}=\frac 1k \sum_{j \sim i} W^n_j
$$
where $(\widetilde{W_1^n},\dots,\widetilde{W_n^n})$ fails to be exchangeable. However, our approach relies on preliminary results like
Theorem \ref{th1} that are general enough to allow us to handle this situation under some mild additional hypothesis. Indeed, assume that
the observations $((x^n_i)_{1 \leq i \leq n})_{n \geq 1}$ satisfy
\begin{equation}\label{eq:xk}
\frac 1n \sum_{i=1}^{n}\delta_{(x_{i}^n,\ldots,x_{i+k-1}^n)} \cvloi \mu^{(k)} \in M_{1}(\Sigma^k).
\end{equation}
\noindent
the $i$ indices being taken modulo $n$. This situation arises e.g. when we are given the realization
$x^n_1,\dots,x^n_n$ of $X^n_1,\dots,X^n_n$ the first $n$ components of 
a stationary Markov chains $(Y_i)_{i \geq 1}$ with transition probability $P$ and stationary measure $\mu$
as in \cite{d76} (see also \cite{b87}). In this case we get
$$
\frac 1n \sum_{i=1}^{n} \delta_{(X^n_i,\dots,X^n_{i+k-1})} \cvloi \mu \otimes \underbrace {P \otimes \cdots  \otimes P}_{k-1} \ \ a.s.
$$
so (\ref{eq:xk}) is satisfied with $\mu^{(k)} = \mu \otimes \underbrace {P \otimes \cdots  \otimes P}_{k-1}$.
which leads to the following conditional LDP 
\begin{theorem}\label{th:kk}
The sequence $(\widetilde{\mathcal L}_n)$ satisfies an LDP on $M_1(\Sigma)$ endowed the weak convergence topology 
with good rate function
$$
\widetilde{\mathcal K}(\nu;\mu^{(k)})=\inf\Big\{\frac1kH(\nu^{(k)}|\mu^{(k)}) : \nu^{(k)}\in M_1(\Sigma^k), 
\frac1k \sum_{i=1}^k \nu_i^{(k)}=\nu \Big\}
$$
\end{theorem}

\noindent
In the particular case of iid observations we further obtain

\begin{corollary} \label{nasri}
The sequence $(\widetilde{\mathcal L}_n)$ satisfies an LDP on $M_1(\Sigma)$ endowed the weak convergence topology 
with good rate function
$$
\widetilde{\mathcal K}(\nu;\mu^{(k)})=H(\nu|\mu).
$$
\end{corollary}

\noindent
hence the $k$-blocks bootstrap is conditionnally efficient in this case 
but fails to be unconditionally efficient for the same reason as Efron's bootstrap, see
Section \ref{fronfron}.


\section{Proof of Theorem \ref{th1}}


Most of the proof of Theorem \ref{th1} relies on the proof of a particular case that we describe now: We are given 
two triangular arrays $((w_{i}^{n})_{1 \leq i \leq n})_{n \geq 1}$
and $((x_{i}^{n})_{1 \leq i \leq n})_{n \geq 1}$ of elements of
$\mathbb{R}_+$ and $\Sigma$ respectively, possibly with repetition, such that 

$$
\nu^{1,n} = \frac{1}{n} \sum_{i=1}^{n} \delta_{w_{i}^{n}} \cvw1 \nu^1 \in M^{1}_{1}(\mathbb{R}_+) \ \ \hbox{and} \ \
\nu^{2,n} = \frac{1}{n} \sum_{i=1}^{n} \delta_{x_{i}^{n}} \cvloi \nu^2 \in M_{1}(\Sigma),
$$
and such that for every $n \geq 1$ we have $\sum_{i=1}^n w^n_i =n$. Let $(\sigma_n)_{n \geq 1}$ be a sequence of 
random variables defined on $(\Omega,\mathcal{A},\mathbb{P})$
such that for every $n \geq 1$ the distribution of $\sigma_n$
is uniform over $\mathfrak{S}_n$.
Let $(T^n)_{n \geq 1}$ be the sequence of random measures defined on $(\Omega,\mathcal{A},\mathbb{P})$ by

$$
T^n = \frac{1}{n} \sum_{i=1}^{n} \delta_{(w_{\sigma_n(i)}^n,x_i^n)} \in \mathcal{M}^1_1(\mathbb{R}_+ \times \Sigma).
$$
This is a particular case of $\mathcal{V}^n$. Following the proof 
of Theorem 1 in \cite{jt08} we shall prove 

\begin{lemma} \label{lemma1} The sequence $(T^n)_{n \geq 1}$ satisfies an LDP on $M_1(\mathbb{R}_+ \times \Sigma)$ 
endowed with the weak convergence topology with good rate function

$$
I(\rho) = \left\{
\begin{array}{cl}
H(\rho| \nu^1 \otimes \nu^2) & \mbox{if \ } \rho_1 = \nu^1 \mbox{and \ } \rho_2 = \nu^2  \\
+ \infty & \mbox{otherwise}.
\end{array}
\right.
$$
\end{lemma}


\noindent
A stronger version of Lemma \ref{lemma1} i.e. 
an LDP for $(T^n)_{n \geq 1}$
on $\mathcal{M}_{1}^{1}(\mathbb{R}_+ \times \Sigma)$ endowed with the distance $\Delta$ is proved in 
Section \ref{stronger}. 
The final proof
of Theorem 1 relies on the latter result and Theorem 2.3 in \cite{g96} and is given in Section \ref{conclusionth1}. 

\bigskip


\subsection{Proof of Lemma \ref{lemma1}}


Let us introduce some more notations before we carry on with the proof of Lemma \ref{lemma1}.
For every $n \geq 1$ we write
$$
\mathcal{P}_n = \left\{ \nu \in \mathcal{M}^1_1(\mathbb{R}_+ \times \Sigma) : \exists \ \sigma \in \mathfrak{S}_n,  \ 
\nu = \frac{1}{n} \sum_{i=1}^{n} \delta_{(w_{\sigma(i)}^{n},x_{i}^{n})} \right\}.
$$
We introduce two triangular arrays $((L_{i}^ {n})_{1 \leq i \leq n})_{n \geq 1}$ 
and $((R_{i}^ {n})_{1 \leq i \leq n})_{n \geq 1}$ of elements of $\mathbb{R}_+ $ and $\Sigma$
respectively, defined on $(\Omega, \mathcal{A}, \P)$ and such that for every $n \geq 1$ the $2n$ random variables 
$L^{n}_{1},\dots, L^{n}_{n}, R^{n}_{1},\dots, R^{n}_{n}$ are mutually independent.
We further assume that every $L^n_i$ (resp. $R^n_i$) is distributed according to $\nu^{1,n}$ 
(resp. $\nu^{2,n}$). The sequence 

$$
\mathcal{T}^n = \frac{1}{n} \sum_{i=1}^{n} \delta_{(L_{i}^{n}, R_{i}^{n})} \in M_1(\mathbb{R}_+ \times \Sigma)
$$
has the following LD behavior


\begin{lemma} \label{lemma2}
The sequence $(\mathcal{T}^n)_{n \geq 1}$ satisfies an LDP on $M_1(\mathbb{R}_+ \times \Sigma)$ endowed 
with the weak convergence 
topology with good rate function $H(\rho | \nu^1  \otimes \nu^2)$.
\end{lemma}


\noindent
{\bf Proof}
Since $\nu^{1,n} \cvw1 \nu^1 $ we have $\nu^{1,n} \cvloi \nu^1 $ according to e.g. Theorem 7.11 in \cite{v03}.
Moreover, since $\nu^{2,n} \cvloi \nu^2$ we have 
$\nu^{1,n} \otimes \nu^{2,n} \cvloi \nu^1  \otimes \nu^2$ 
(see \cite{b68}, Chapter 1, Theorem 3.2). The announced result then follows from Theorem 3 
in \cite{bj88}. \hfill $\Box$

\bigskip

\noindent
Our strategy in proving Lemma \ref{lemma1} consists in comparing $T^n$ to random measures 
coupled to $\mathcal{T}^n$. Comparison 
is possible because the $\rho \in M_1(\mathbb{R}_+ \times \Sigma)$ such that $I(\rho) < + \infty$ can 
be approached in the weak convergence topology by elements of $\mathcal{P}_n$. Our proof of this property
requires to use several metrics on $M_1(\mathbb{R}_+ \times \Sigma)$ compatible with the weak convergence topology. 
This is the reason why in Section \ref{approche} we give a short account on the weak convergence topology 
prior to the proof of our approximation result. In Section \ref{couplage}
we construct our coupling. Finally, in Section \ref{ineg} we prove the LD bounds of Lemma \ref{lemma1}.


\subsubsection{An approximation result} \label{approche}


We are given the Polish space $(\mathbb{R}_+ ,\dg)$ with $d_e(x,y)=|x-y|$\footnote{We use this notation
to underscore the fact that for Lemma \ref{lemma1} to hold it is not necessary for the $w^n_i$'s to be real numbers. 
It is sufficient to have Polish space valued $w^n_i$'s.}. The distance $\dg$ is not a 
bounded metric but it is topologically equivalent to the bounded metric

$$
\tilde{\dg}(x,y) = \frac{\dg(x,y)}{1 + \dg(x,y)}
$$
Analogously
we define a bounded metric $\tilde{\de}$ on $\Sigma$ on the ground of $\de$.
The product topology on $\mathbb{R}_+ \times \Sigma$ is metrizable by e.g. 

\begin{equation}
d_{2,M}((x_1,x_2),(y_1,y_2)) = \max ( \dg (x_1,y_1), \de (x_2,y_2))
\label{dm}
\end{equation}
or 

\begin{equation}
d_{2,+}((x_1,x_2),(y_1,y_2)) = \dg (x_1,y_1) + \de (x_2,y_2).
\label{d+}
\end{equation}
They both make $\mathbb{R}_+ \times \Sigma$ a Polish space. We can also metrize 
the product topology on $\mathbb{R}_+ \times \Sigma$ with the analogues 
$\tilde{d}_{2,M}$ and $\tilde{d}_{2,+}$ of (\ref{dm}) and (\ref{d+}) built on the ground of $\tilde{\dg}$
and $\tilde{\de}$. With a slight abuse of notation we shall denote by

\begin{equation}
\beta_{BL,\delta}(\rho,\nu) = 
\sup_{f \in  C_{b}(\mathbb{R}_+ ) \atop \|f\|_{\infty} + \|f\|_{L,\delta} \leq 1} 
\left\{ \left| \int_{\mathbb{R}_+ } f d\rho - \int_{\mathbb{R}_+ } f d\nu \right| \right\}  \label{beta}
\end{equation}

\noindent
the so-called dual-bounded Lipschitz 
metric on $M_1(\mathbb{R}_+)$ where $\delta$ is either $d_e$ or $\tilde{d}_e$. It is compatible with the weak convergence topology 
(see \cite{d02}, Chapter 11, Theorem 11.3.3).
According to Kantorovitch-Rubinstein Theorem (see \cite{d02}, Chapter 11, Theorem 11.8.2)
the following metric on $M_1(\mathbb{R}_+ )$

$$
W_{1,\tilde{\dg}} (\rho,\nu) = \inf_{Q \in \mathcal{C}(\rho,\nu)} 
\left\{ \int_{\mathbb{R}_+ \times \mathbb{R}_+} \tilde{d}_{e}(x,y) Q(dx,dy)   \right\},
$$
the so-called Wasserstein-1 metric associated to $\tilde{\dg}$ is compatible with the weak convergence topology as well.
However, note that the "analogue" of $W_{1,\tilde{\dg}}$ built on the ground of $\dg$ is not a 
metric for the weak convergence 
topology (for an illustration of this fact see \cite{d02} p.420-421). Finally we shall denote by 

\begin{equation}
\beta_{BL,\tilde{d}_{2,M}}(\rho,\nu) = 
\sup_{f \in  C_{b}(\mathbb{R}_+ \times \Sigma) \atop \|f\|_{\infty} + \|f\|_{L,\tilde{d}_{2,M}} \leq 1} 
\left\{ \left| \int_{\mathbb{R}_+ \times \Sigma} f d\rho - 
\int_{\mathbb{R}_+ \times \Sigma} f d\nu \right| \right\}  \label{betap}
\end{equation}

\noindent
and  

$$
W_{1,\tilde{d}_{2,+}}(\rho,\nu) = \inf_{Q \in \mathcal{C}(\rho,\nu)} 
\left\{ \int_{(\mathbb{R}_+ \times \Sigma) \times (\mathbb{R}_+ \times \Sigma)} \tilde{d}_{2,+}(x,y) Q(dx,dy) \right\},
$$
two metrics on $M_1(\mathbb{R}_+ \times \Sigma)$ compatible with the weak convergence 
topology on this set. The following is a key result in the proof of Lemma \ref{lemma1}. 


\begin{lemma}
Let $\rho \in M_1(\mathbb{R}_+ \times \Sigma)$ be such that $\rho_1 = \nu^1 $ and $\rho_2 = \nu^2$. 
For every $n \geq 1$ there exists a
$\rho_n \in \mathcal{P}_n$ such that $\rho_n \cvloi \rho$.
\label{lem2}
\end{lemma}


\noindent
{\bf Proof}
Let $\rho \in M_1(\mathbb{R}_+ \times \Sigma)$ be such that $\rho_1 = \nu^1 $ and $\rho_2 = \nu^2$. According to 
Varadarajan's Lemma (see \cite{d02}, Chapter 11,
Theorem 11.4.1) there exists a family $((u_{i}^{n},v_{i}^{n} )_{1 \leq i \leq n})_{n \geq 1}$ of 
elements of $\mathbb{R}_+ \times \Sigma$ such that

$$
\gamma^n  = \frac{1}{n} \sum_{i=1}^{n} \delta_{(u_{i}^{n}, v_{i}^{n})} \cvloi \rho.
$$
For every $n \geq 1$ we take $\varphi_n, \tau_n \in \mathfrak{S}_n$ such that

$$
\sum_{i=1}^{n} \tilde{d_e}(u_{i}^{n}, w_{\varphi_n(i)}^{n})
= \min_{\varphi \in \mathfrak{S}_n} \left\{  \sum_{i=1}^{n} 
\tilde{d_e}(u_{i}^{n}, w_{\sigma(i)}^{n}) \right\}
$$

\noindent
and

$$
\sum_{i=1}^{n} \tilde{d_{\Sigma}}(v_{i}^{n}, x_{\tau_n(i)}^{n})
= \min_{\tau \in \mathfrak{S}_n} \left\{  \sum_{i=1}^{n} 
\tilde{d_{\Sigma}}(v_{i}^{n}, x_{\tau(i)}^{n}) \right\}.
$$
We shall prove that 

$$
\rho^n = \frac{1}{n} \sum_{i=1}^{n} \delta_{(w_{\varphi_n(i)}^{n},x_{\tau_n(i)}^{n})} \cvloi \rho.
$$
To this end it is sufficient to prove that $\beta_{BL,\tilde{d}_{2,M}}(\rho^n,\rho) \rightarrow 0$. 
Let $\varepsilon > 0$ be fixed. There exists an $N_0$ such that for every $n \geq N_0$ we have

\begin{equation}	
\beta_{BL, \tilde{d}_{2,M}}(\rho, \gamma^n ) <  \varepsilon/3.
\label{pe}
\end{equation}

\noindent
Since $\gamma^n_1 \cvloi \rho_1 = \nu^1 $ there exists an $N_1$ such that for 
every $n \geq N_1$ we have $W_{1,\tilde{\dg}}(\gamma^n_1, \nu^{1,n}) < \varepsilon / 3$.
We will show that due to this for every $n \geq N_1$ we have

\begin{equation}
\frac{1}{n} \sum_{i=1}^{n} \tilde{d_e}(u_{i}^{n}, w_{\varphi_n(i)}^{n}) < \varepsilon / 3.
\label{se}
\end{equation}
We shall prove that it leads to 

\begin{equation}
\beta_{BL, \tilde{d}_{2,M}}(\hat{\gamma}^n, \gamma^n) < \varepsilon / 3
\label{te}
\end{equation}
for every $n \geq N_1$ where

$$\hat{\gamma}^n = \frac{1}{n} \sum_{i=1}^{n} 
\delta_{(w_{\varphi_n(i)}^{n}, v_{i}^{n})}.$$ 
Analogously  to (\ref{se} - \ref{te}) one can prove that there exists an $N_2$ such that for every $n \geq N_2$

\begin{equation}
\beta_{BL,\tilde{d}_{2,M}}( \hat{\gamma}^n, \rho^n) 
< \varepsilon / 3.
\label{qe}
\end{equation}
By combining (\ref{pe}, \ref{te}, \ref{qe}) we obtain the announced result. \\
\bigskip
\underline{{\it Proof of (\ref{se})}}
Since $\gamma_1^n$ and $\nu^{1,n}$ have finite support every Borel probability 
measure $Q$ on $\mathbb{R}_+ \times \mathbb{R}_+$ such that $Q_1 = \gamma_{1}^{n}$ and $Q_2 = \nu^{1,n}$ is of the form

\begin{equation}
Q(\alpha) = \frac{1}{n} \sum_{i,j=1}^{n} \alpha_{i,j} \delta_{(u_{i}^{n},w_{j}^{n})}
\label{bisto}
\end{equation}
where $\alpha = (\alpha_{i,j})_{1 \leq i,j \leq n}$ is an $n \times n$ bi-stochastic matrix.
Conversely every $n \times n$ bi-stochastic matrix $\alpha$ defines through (\ref{bisto}) 
a Borel probability measure $Q(\alpha)$ on $\mathbb{R}_+ \times \mathbb{R}_+$ such that $Q(\alpha)_1 = \gamma_{1}^{n}$ 
and $Q(\alpha)_2 = \nu^{1,n}$. From the Birkhoff-Von Neumann Theorem 
we know that every bi-stochastic matrix can be written
as a convex combination of permutation matrices. These are $n \times n$ matrices 
with a single 1 in every line and every column, all other entries being 0 (for a proof of this 
fact see e.g. \cite{l68}, Chapter 11, Example 11.2). There is an obvious one-to-one correspondence 
between elements of $\mathfrak{S}_n$ and $n \times n$ permutation 
matrices. For every $\delta \in \mathfrak{S}_n$ we shall denote by 
$K_{\delta}$ the permutation matrix naturally associated to $\delta$ by this correspondence. Therefore
for any Borel probability measure $Q$ on $\mathbb{R}_+ \times \mathbb{R}_+ $ such that
$Q_1 = \gamma_{1}^{n}$ and $Q_2 = \nu^{1,n}$
i.e. any choice of the components $(\lambda_{\delta})_{\delta \in \mathfrak{S}_n}$ of the convex 
combination $\alpha = \sum_{\delta \in \mathfrak{S}_n} \lambda_{\delta} K_{\delta}$ such that 
$Q = Q(\alpha)$ we have

\begin{eqnarray}
\int_{\mathbb{R}_+ \times \mathbb{R}_+} \tilde{d}_e(x,y) Q(dx,dy) & = & \sum_{\delta \in \mathfrak{S}_n} 
                                              \lambda_{\delta} 
                                              (\sum_{i=1}^{n} \tilde{d}_e(u_{i}^{n}, w_{\delta(i)}^{n})) \\
                                        & \geq & \frac{1}{n} \sum_{i=1}^{n} 
                                        \tilde{d}_e(u_{i}^{n}, w_{\varphi_n(i)}^{n}).
\end{eqnarray}
Hence for every $n \geq N_1$

\begin{eqnarray*}
\frac{1}{n} \sum_{i=1}^{n} \tilde{d}_{e}(u_{i}^{n}, w_{\varphi_n(i)}^{n}) & = &
W_{1,\tilde{\dg}}( \gamma_1^n, \nu^{1,n})
 \\
& < & \varepsilon / 3
\end{eqnarray*}

\noindent
which proves (\ref{se}). \\
\bigskip
\underline{{\it Proof of (\ref{te})}} For every $f \in \Cb$ such that 
$||f||_{\infty} + ||f||_{L,\tilde{d}_{2,M}} \leq 1$ we have

\begin{eqnarray*} 
\left| \int_{\mathbb{R}_+ \times \Sigma} f d \hat{\gamma }^n - \int_{\mathbb{R}_+ \times \Sigma} f d \gamma^n \right| 
& \leq & \frac{1}{n} \sum_{i=1}^{n}
\left|f(w_{\varphi_n(i)}^{n}, v_{i}^{n}) - f(u_{i}^{n}, v_{i}^{n}) \right| \\
& \leq & \frac{1}{n} \sum_{i=1}^{n} \tilde{d}_{2,M} ((w_{\varphi_n(i)}^{n}, v_{i}^{n}), (u_{i}^{n}, v_{i}^{n})) \\ 
& \leq & \frac{1}{n} \sum_{i=1}^{n} \tilde{\dg}(u_{i}^{n}, w_{\varphi_n(i)}^{n}) \\
& < & \varepsilon / 3.
\end{eqnarray*}

\noindent
Hence for every $n \geq N_1$ we have 
$\beta_{BL, \tilde{d}_{2,M}}(\hat{\gamma}^n, \gamma^n) < \varepsilon / 3$.  \hfill $\Box$ 


\subsubsection{Coupled empirical measures} \label{couplage}


To every $n \geq 1$ and every realization of $\mathcal{T}^n$ we associate two 
elements of $M_1(\mathbb{R}_+ \times \Sigma)$ by 

\begin{equation} 
W_{1,\tilde{d}_{2,+}}(\mathcal{T}^n,\widetilde{W}^n) = 
\min_{\nu \in \mathcal{P}_n} \left\{ W_{1,\tilde{d}_{2,+}}(\mathcal{T}^n,\nu) \right\}
\label{argmin}
\end{equation}

\noindent
and  

\begin{equation} 
W_{1,\tilde{d}_{2,+}}(\mathcal{T}^n,\widehat{W}^n) = 
\max_{\nu \in \mathcal{P}_n} \left\{ W_{1,\tilde{d}_{2,+}}(\mathcal{T}^n,\nu) \right\}.
\label{argmax}
\end{equation}

\noindent
In case there are several elements of $\mathcal{P}_n$ achieving the 
min (resp. the max) $\widetilde{W}^n$ (resp. $\widehat{W}^n$) is 
picked uniformly at random among these measures.

\begin{lemma} For every $n \geq 1$ the random measures $\widetilde{W}^n, \widehat{W}^n$ and $T^n$ are 
identically distributed over $M_1(\mathbb{R}_+ \times \Sigma)$. 
\label{lem3}
\end{lemma}

\noindent
{\bf Proof}
We shall only prove that $\widetilde{W}^n$ and $T^n$ are identically distributed 
since the proof with $\widehat{W}^n$ and $T^n$
works the same way. Let $n \geq 1$ be fixed. For the sake of clarity let us first assume that 
there is no repetition among the $w_1^n,\dots,w_n^n$ and the 
$x_1^n,\dots,x_n^n$. In this case every $\rho \in \mathcal{P}_n$ corresponds to a single $\tau \in \mathfrak{S}_n$ by

\begin{equation}
\rho = \frac{1}{n} \sum_{i=1}^{n} \delta_{(w_{\tau(i)}^{n},x_{i}^{n})}.
\label{code}
\end{equation}
Next let us consider a fixed realization $(l_i^n,r_i^n)_{1 \leq i \leq n}$ of $(L_i^n,R_i^n)_{1 \leq i \leq n}$
and let us denote by w$^n$ the corresponding value for $\mathcal{T}^n$. Due to the property of the minimizer in the 
Wasserstein distance between two atomic measures we already employed in the proof of Lemma \ref{lem2}, for 
every $\rho \in \mathcal{P}_n$ (i.e. every $\tau \in \mathfrak{S}_n$ according to (\ref{code})) there 
exists a $\sigma \in \mathfrak{S}_n$ such that
\begin{eqnarray*}
W_{1,\tilde{d}_{2,+}}({\mbox w}^n,\rho) & = & \frac{1}{n} \sum_{i=1}^{n} \tilde{d}_{2,+} ((l_i^n,r_i^n), 
(w_{\sigma(i)}^{n}, x_{\sigma \circ \tau(i)}^{n})) \\
 & = & \frac{1}{n} \sum_{i=1}^{n} \tilde{d}_{e}(l_i^n,w_{\sigma(i)}^{n}) 
+ \frac{1}{n} \sum_{i=1}^{n} \tilde{d_{\Sigma}}(r_i^n,x_{\sigma \circ \tau(i)}^{n}).
\end{eqnarray*}

\noindent
Thus, since we are looking for the minimum over $\sigma$ and $\tau$, for a fixed 
realization $(l_i^n,r_i^n)_{1 \leq i \leq n}$ of $(L_i^n,R_i^n)_{1 \leq i \leq n}$ the 
corresponding value of $\widetilde{{\mbox w}}^n$ is obtained by finding $\eta_l$ and $\eta_r$ such that 

$$
\sum_{i=1}^{n} \tilde{d}_e (l_i^n,w_{\eta_l(i)}^n) = \min_{\sigma \in \mathfrak{S}_n} \left\{ \sum_{i=1}^{n} 
\tilde{d}_e (l_i^n,w_{\sigma(i)}^n) \right\}
$$
and
$$
\sum_{i=1}^{n} \tilde{d_{\Sigma}} (r_i^n,x_{\eta_r(i)}^n) = \min_{\sigma \in \mathfrak{S}_n} \left\{ \sum_{i=1}^{n} 
\tilde{d_{\Sigma}} (r_i^n,x_{\sigma(i)}^n) \right\}
$$
and taking

$$
\widetilde{{\mbox w}}^n = \frac{1}{n} \sum_{i=1}^{n} \delta_{(w_{\eta_l(i)}^n,x_{\eta_r(i)}^n)}.
$$
In case several $\eta_l$ and/or $\eta_r$ realize the minima in the displays above, those defining 
$\widetilde{{\mbox w}}^n$ 
are picked among them uniformly at random. Now, remark that for every $\gamma_{l},\gamma_{r} \in \mathfrak{S}_n$,  
observing $(l_{\gamma_l(i)}^{n}, r_{\gamma_r(i)}^{n})_{1 \leq i \leq n}$ has the same probability as observing  
$(l_i^n,r_i^n)_{1 \leq i \leq n}$ and results in $\gamma_l \circ \eta_l$ and $\gamma_r \circ \eta_r$ in defining
$\widetilde{{\mbox w}}^n$ instead of $\eta_l$ and $\eta_r$. Thus, if 
we consider $\eta_l$ and $\eta_r$ as random variables 
defining $\widetilde{W}^n$, we see that their distribution {\it conditioned on }$W^n$ is uniform over $\mathfrak{S}_n$.
Hence $\widetilde{W}^n$ and $T^n$ are both uniformly 
distributed over $\mathcal{P}^n$, thus identically distributed. 
This proof extends easily to the case when there are repetitions among the $w_1^n,\dots,w_n^n$
or $x_1^n,\dots,x_n^n$. \hfill $\Box$


\subsubsection{Proof of the LD bounds of Lemma 1} \label{ineg}


We start the proof of the LD bounds by proving the following

\begin{lemma} We have:
\begin{enumerate}
\item $I$ is a good rate function.
\item The sequence $(T^n)_{n \geq 1}$ is exponentially tight. 
\end{enumerate}
\label{lem4}
\end{lemma}

\noindent
{\bf Proof}\\
(1) Let $\alpha \geq 0$. We have

\begin{eqnarray*}
N_{\alpha} & = & \{ \rho \in M_1(\mathbb{R}_+ \times \Sigma) : I(\rho) \leq \alpha \} \\
& = & \{ \rho \in M_1(\mathbb{R}_+ \times \Sigma) : H(\rho | \nu^1  \otimes \nu^2) \leq \alpha \} \cap 
\{ \rho \in M_1(\mathbb{R}_+ \times \Sigma) : \rho_1 = \nu^1  \mbox{ and } \rho_2 = \nu^2 \}.
\end{eqnarray*}

\noindent
Thus, for every $\alpha \geq 0$, $N_{\alpha}$ is the intersection of a compact and a closed subset of 
$M_1(\mathbb{R}_+ \times \Sigma)$, therefore it is compact.\\

\noindent
(2) For every measurable $A \subset M_1(\mathbb{R}_+ \times \Sigma)$ we have 

\begin{eqnarray*}
\limsup_{n \rightarrow \infty} \frac{1}{n} \log \mathbb{P}( T^n \in A^c) & = & \limsup_{n \rightarrow \infty} 
\frac{1}{n} \log \P( \mathcal{T}^n \in A^c | \frac{1}{n} \sum_{i=1}^{n} 
\delta_{L_i^n} = \nu^{1,n} , \frac{1}{n} \sum_{i=1}^{n} \delta_{R_i^n} = \nu^{2,n} ) \\
& \leq & \limsup_{n \rightarrow \infty} \frac{1}{n} \log \P( \mathcal{T}^n \in A^c) \\ 
& & \ \ \ \ \ \ \ \ -  \liminf_{n \rightarrow \infty} \frac{1}{n} 
\log \P(\frac{1}{n} \sum_{i=1}^{n} \delta_{L_i^n} = \nu^{1,n} ) \\ 
& & \ \ \ \ \ \ \ \ -  \liminf_{n \rightarrow \infty} \frac{1}{n} 
\log \P(\frac{1}{n} \sum_{i=1}^{n} \delta_{R_i^n} = \nu^{2,n} ).
\end{eqnarray*}

\noindent
Since $(\mathcal{T}^n)_{n \geq 1}$ satisfies an LDP on $M_1(\mathbb{R}_+ \times \Sigma)$ with a good 
rate function it is 
exponentially tight (see \cite{dz98}, Remark a) p.8). Thus for every $\alpha \geq 0$ we can chose a 
compact set $A_{\alpha} \subset M_1(\mathbb{R}_+ \times \Sigma)$ that makes the first term in the last display 
smaller than $ - \alpha - 2$. On the other hand it is clear that
for every $n \geq 1$ we have

$$
\P( \frac{1}{n} \sum_{i=1}^{n} \delta_{L_i^n} = \nu^{1,n} ) \geq n!  \frac{1}{n^n}
$$
equality corresponding to the case when there are no ties among the $x^n_1,\dots, x^n_n$. Thus

$$
-  \liminf_{n \rightarrow \infty} \frac{1}{n} \log \P(\frac{1}{n} \sum_{i=1}^{n} \delta_{L_i^n} = \nu^{1,n} ) \leq 1
$$
which completes the proof. \hfill $\Box$

\bigskip

\noindent
\emph{Proof of the lower bound}

\noindent
It is sufficient in order to prove the lower bound of the LDP to prove that 

$$
- I(\rho) \leq \liminf_{n \rightarrow \infty} \frac{1}{n} \log \mathbb{P}( T^n \in B(\rho,\varepsilon))
$$
holds for every $\rho \in M_1(\mathbb{R}_+ \times \Sigma)$ and every $\varepsilon > 0$, 
where $B(\rho,\varepsilon)$ stands for 
the open ball centered 
at $\rho \in M_1(\mathbb{R}_+ \times \Sigma)$ of radius $\varepsilon > 0$ for the $W_{1,\tilde{d}_{2,+}}$ 
metric. So let $\varepsilon > 0$ and 
$\rho \in M_1(\mathbb{R}_+ \times \Sigma)$ be such that $I(\rho) < + \infty$. In particular $\rho_1 = \nu^1 $ and
$\rho_2 = \nu^2$. 
According to Lemma \ref{lem2} there exists a sequence $(\rho^n)_{n \geq 1}$ of elements 
of $M_1(\mathbb{R}_+ \times \Sigma)$ 
such that $\rho^n \in \mathcal{P}_n$ and $\rho^n \cvloi \rho$. According to Lemma \ref{lem3}

\begin{eqnarray*}
\mathbb{P}(T^n \in B(\rho,\varepsilon)) & = & \P( \widetilde{W}^n \in B(\rho,\varepsilon)) \\
& \geq & \P( W_{1,\tilde{d}_{2,+}}(\widetilde{W^n},\mathcal{T}^n) < \frac{\varepsilon}{3},  
W_{1,\tilde{d}_{2,+}}(\mathcal{T}^n,\rho^n) < \frac{\varepsilon}{3}, 
W_{1,\tilde{d}_{2,+}}(\rho^n,\rho) < \frac{\varepsilon}{6}) \\
& \geq & \P( W_{1,\tilde{d}_{2,+}}(\rho^n,\mathcal{T}^n) < \frac{\varepsilon}{3}, 
W_{1,\tilde{d}_{2,+}}(\rho^n,\rho) < \frac{\varepsilon}{6}) \\
\end{eqnarray*}

\noindent
since it follows from the definition of $\widetilde{W^n}$ that for every $\rho^n \in \mathcal{V}_n$ we have

$$
W_{1,\tilde{d}_{2,+}}(\widetilde{W^n},\mathcal{T}^n) \leq W_{1,\tilde{d}_{2,+}}(\rho^n,\mathcal{T}^n).
$$
On the other hand since $\rho^n \cvloi \rho$ we get that for $n$ large enough 
$\left\{ W_{1,\tilde{d}_{2,+}}(\rho^n,\rho) < \frac{\varepsilon}{6} \right\} = \Omega $. Thus, for those $n$'s

\begin{eqnarray*}
\P( W_{1,\tilde{d}_{2,+}}(\rho^n,\mathcal{T}^n) < \frac{\varepsilon}{3}, 
W_{1,\tilde{d}_{2,+}}(\rho^n,\rho) < \frac{\varepsilon}{6}) & \geq & 
\P( W_{1,\tilde{d}_{2,+}}(\rho,\mathcal{T}^n) < \frac{\varepsilon}{6}, 
W_{1,\tilde{d}_{2,+}}(\rho^n,\rho) < \frac{\varepsilon}{6}) \\
& \geq & \P( W_{1,\tilde{d}_{2,+}}(\mathcal{T}^n,\rho) < \frac{\varepsilon}{6}).
\end{eqnarray*}

\noindent
Finally, it follows from Lemma \ref{lemma1} that
\begin{eqnarray*}
\liminf_{n \rightarrow \infty} \frac{1}{n} \log \mathbb{P}( T^n \in B(\rho,\varepsilon)) & \geq & 
\liminf_{n \rightarrow \infty} \frac{1}{n} 
\log \P( W_{1,\tilde{d}_{2,+}}(\rho,\mathcal{T}^n) 
< \frac{ \varepsilon}{6}) \\
& \geq & - H(\rho | \nu^1  \otimes \nu^2) = - I(\nu).
\end{eqnarray*}

\noindent
\emph{Proof of the upper bound}

\noindent
In order to prove the upper bound of the LDP, it is sufficient to prove that it holds for compact subsets of
$M_1(\mathbb{R}_+ \times \Sigma)$. Indeed, since $(T^n)_{n \geq 1}$ is an exponentially tight 
sequence (see Lemma \ref{lem4}) 
the full upper bound will follow from Lemma 1.2.18 in \cite{dz98}. Let $A$ be a compact subset 
of $M_1(\mathbb{R}_+ \times \Sigma)$ and let us denote by

$$
A_{\nu^1 , \nu^2} = \left\{ \rho \in A  : \rho_1 = \nu^1  \mbox{ and } \rho_2 = \nu^2 \right\}
$$
which is a compact subset of $M_1(\mathbb{R}_+ \times \Sigma)$ as well. Since the weak convergence 
topology on $M_1(\mathbb{R}_+ \times \Sigma)$ is compatible with the $W_{1,\tilde{d}_{2,+}}$ metric, it 
makes $M_1(\mathbb{R}_+ \times \Sigma)$ 
a regular topological space: For every $\rho \in A$ such 
that $\rho \in A_{\nu^1 , \nu^2}^c$ there exists $\varepsilon_{\rho} > 0$ such 
that $B(\rho, 2 \varepsilon_{\rho})  \cap A_{\nu^1 , \nu^2} = \emptyset$. In particular
$\bar{B}(\rho, \varepsilon_{\rho})  \cap A_{\nu^1 , \nu^2} = \emptyset$ where $\bar{B}(\rho, \varepsilon)$ denotes the
closed ball centered on $\rho \in M_1(\mathbb{R}_+ \times \Sigma)$ of radius $\varepsilon > 0$ for 
the $W_{1,\tilde{d}_{2,+}}$ metric. 
On the other hand, since $\rho \mapsto H(\rho | \nu^1  \otimes \nu^2)$ is lower semi-continuous, for 
every $\rho \in A_{\nu^1 , \nu^2}$ and every $\delta > 0$ there exists a $\varphi(\rho, \delta) > 0$ such that

$$
\inf_{\gamma \in \bar{B}(\rho, \varphi(\rho, \delta))} H(\gamma | \nu^1  \otimes \nu^2) \geq (H(\rho | \nu^1  \otimes \nu^2) - \delta) 
\wedge \frac{1}{\delta}.
$$
For every $\delta > 0$ we consider the coverage 

$$
A \subset \left( \cup_{\rho \in A \cap A_{\nu^1 , \nu^2}^c} B (\rho, \varepsilon_{\rho}) \right) \cup \left( \cup_{\rho \in  A_{\nu^1 , \nu^2}} 
B (\rho, \frac{\varphi(\rho, \delta)}{8}) \right) 
$$
from which we extract a finite coverage

$$
A \subset \left( \cup_{\rho \in I_1} B (\rho, \varepsilon_{\rho}) \right) \cup \left( \cup_{\rho \in I_2} 
B (\rho, \frac{\varphi(\rho, \delta)}{8}) \right) 
$$
where $I_1 \subset A \cap A_{\nu^1 , \nu^2}^c$ and $I_2 \subset A_{\nu^1 , \nu^2}$ are finite sets. Then, according to Lemma 1.2.15
in \cite{dz98} 

\begin{eqnarray*} 
\limsup_{n \rightarrow \infty} \frac{1}{n} \log \mathbb{P}( T^n \in A) & \leq & \max 
\left\{ \limsup_{n \rightarrow \infty} \frac{1}{n} \log \mathbb{P}( T^n \in  
\cup_{\rho \in I_1} \bar{B} (\rho, \varepsilon_{\rho}) \cap A), \right. \\
& & \ \ \ \ \left. \limsup_{n \rightarrow \infty} \frac{1}{n} \log \mathbb{P}( T^n \in  \cup_{\rho \in I_2} 
B (\rho, \frac{\varphi(\rho, \delta)}{8})) \right\}.
\end{eqnarray*}

\noindent
For every $\rho \in I_1$ there can not be an infinite number of integers $n_k$ such that
$$\mathbb{P}( T^{n_k} \in  
\cup_{\rho \in I_1} \bar{B} (\rho, \varepsilon_{\rho}) \cap A) \neq 0$$
for otherwise we would get 
$\bar{B} (\rho, \varepsilon_{\rho}) \cap A_{\nu^1 , \nu^2} \neq \emptyset$.
The first term in the max is then equal to $- \infty$. We are left with the second term 
and according to Lemmas 1,2 and 3
we have

\begin{eqnarray*}
\limsup_{n \rightarrow \infty} \frac{1}{n} \log \mathbb{P}( T^n \in A) & \leq & 
\limsup_{n \rightarrow \infty} \frac{1}{n} \log 
\mathbb{P}( T^n \in  \cup_{\rho \in I_2} B (\nu, \frac{\varphi(\rho, \delta)}{8})) \\
& \leq & 
\max_{\rho \in I_2} \left\{ \limsup_{n \rightarrow \infty} \frac{1}{n} \log \P( \widehat{W}^n \in B (\rho, \frac{\varphi(\rho, \delta)}{8})) 
\right\} \\
& \leq & \max_{\rho \in I_2} \left\{ \limsup_{n \rightarrow \infty} \frac{1}{n} \log \P( W_{1,\tilde{d}_{2,+}}(\rho^n,\widehat{W}^n) 
< \frac{\varphi(\rho, \delta)}{4}) \right\} \\
& \leq & \max_{\rho \in I_2} \left\{ \limsup_{n \rightarrow \infty} \frac{1}{n} \log \P( W_{1,\tilde{d}_{2,+}}(\rho, W^n) 
< \frac{\varphi(\rho, \delta)}{2}) \right\} \\
& \leq & \max_{\rho \in I_2} \left\{ - \inf_{\gamma \in \bar{B}(\rho, \varphi(\rho, \delta))} H(\gamma | \nu^1  \otimes \nu^2) \right\}\\
& \leq & \max_{\rho \in I_2} \left\{ - (H(\rho | \nu^1  \otimes \nu^2) - \delta) \wedge \frac{1}{\delta} \right\} \\
& \leq & \max_{\rho \in I_2} \left\{ - (I(\rho) - \delta) \wedge \frac{1}{\delta} \right\} \\
& \leq & - \inf_{\rho \in A} \left\{ (I(\rho) - \delta) \wedge \frac{1}{\delta} \right\}.
\end{eqnarray*}

\noindent
By letting $\delta \rightarrow 0$ we obtain the announced upper bound, see Remark 1.2.10 in \cite{dz98}. 


\subsection{A stronger version of Lemma \ref{lemma1}} \label{stronger}


In the present section we shall prove the following
 
\begin{lemma} \label{lemma2bis} The sequence $(T^n)_{n \geq 1}$ satisfies an LDP 
on $\mathcal{M}_1^1(\mathbb{R}_+ \times \Sigma)$ 
endowed with the distance $\Delta$ with good rate function $I$.
\end{lemma}

\noindent
\textbf{Proof} First let us notice that 

$$\mathcal{N}^{1}_{1}(\mathbb{R}_+ \times \Sigma) = \left\{ \rho \in M_{1}(\mathbb{R}_+ \times \Sigma) \ :
 \ \int_{\mathbb{R}_+} x \ \rho_1(dx) \leq 1  \right\}$$
is a closed subset of 
$M_1(\mathbb{R}_+ \times \Sigma)$ when the latter is endowed with the weak convergence topology.
Since for every $n \geq 1$ we have $\mathbb{P}(T_n \in \mathcal{N}^{1}_{1}(\mathbb{R}_+ \times \Sigma)) = 1$, 
Lemma 4.1.5 in \cite{dz98} implies that $(T_n)_{n \geq 1}$ obeys an LDP 
on $\mathcal{N}^{1}_{1}(\mathbb{R}_+ \times \Sigma)$
endowed with the weak convergence topology, with good rate function $I$. 

\bigskip
\noindent
Next we prove that the same remains true when 
$\mathcal{N}^{1}_{1}(\mathbb{R}_+ \times \Sigma)$ is endowed with the distance $\Delta$. Indeed,
since $(T^n)_{n \geq 1}$ satisfies an LDP on the Polish space $M_1(\mathbb{R}_+ \times \Sigma)$ with a good 
rate function it is exponentially tight: For every $L > 0$ there exists a compact 
$A \subset M^1(\mathbb{R}_+ \times \Sigma)$ such that 

$$
\limsup_{n \rightarrow \infty} \frac{1}{n} \log \mathbb{P}( T^n \in A^c) \leq -L.
$$
The set $A \cap K$ with 
$$
K= \left\{ \rho \in M^1(\mathbb{R}_+ \times \Sigma) \ : \rho_1 \in \cup_{i=1}^n \{\nu^{1,n} \} \right\}
$$
is a compact subset of $\mathcal{N}^{1}_{1}(\mathbb{R}_+ \times \Sigma)$ endowed with $\Delta$.
Moreover

\begin{eqnarray*}
\limsup_{n \rightarrow \infty} \frac{1}{n} \log \mathbb{P}( T^n \in (A \cap K)^c) & \leq & 
\limsup_{n \rightarrow \infty} \frac{1}{n} \log \mathbb{P}( T^n \in A^c) + \\
& & \ \ \ \ \ \ \ \ + \limsup_{n \rightarrow \infty} \frac{1}{n} \log \mathbb{P}( T^n \in K^c) \\
 & \leq & -L.
\end{eqnarray*}
 
\bigskip

\noindent
Finally, $\mathcal{M}^{1}_{1}(\mathbb{R}_+ \times \Sigma)$ 
is a closed subset of $(\mathcal{N}^{1}_{1}(\mathbb{R}_+ \times \Sigma),\Delta)$ and for every $n \geq 1$ we have
$\mathbb{P}(T_n \in \mathcal{M}^{1}_{1}(\mathbb{R}_+ \times \Sigma)) = 1$ so again, due to Lemma 4.1.5 in \cite{dz98}
the sequence $(T_n)_{n \geq 1}$ obeys an LDP on $\mathcal{M}^{1}_{1}(\mathbb{R}_+ \times \Sigma)$
endowed with the distance $\Delta$ with good rate function $I$. \hfill $\Box$


\subsection{Conclusion of the proof of Theorem 1} \label{conclusionth1}


In order to conclude the proof of Theorem 1 it is sufficient to establish that the distribution of $\mathcal{V}^n$
on $\mathcal{M}^{1}_{1}(\mathbb{R}_+ \times \Sigma)$ is a mixture of Large Deviation Systems (LDS) in the sense of 
\cite{g96}. For the sake of clarity we recover the notations of \cite{g96} when identifying the components of the LDS

\noindent
$\bullet$ $\mathcal{Z} = \mathcal{M}^{1}_{1}(\mathbb{R}_+ \times \Sigma)$ is a Polish space when 
endowed with the distance
$\Delta$, see the Appendix.\\
$\bullet$ $\mathcal{X} = M^{1}_{1}(\mathbb{R}_+)$ is a Polish space when endowed with the $W_1$-distance since it is a 
closed subset of the Polish space $(\mathcal{W}^{1}(\mathbb{R}_+),W_{1,|\cdot|})$ (see e.g. \cite{bo08}).\\
$\bullet$ For every $n \geq 1$ we note
$$\mathcal{X}_n = \left\{ \nu \in M^{1}_{1}(\mathbb{R}_+) \ : \ \exists \ (w_1,\dots,w_n) \in (\mathbb{R}_+)^n, \ 
\nu = \frac{1}{n} \sum_{i=1}^n \delta_{w_i} \right\}$$
and  for every $\nu \in \mathcal{X}$ and every $n \geq 1$ there exists a $\nu^n \in \mathcal{X}_n$
such that $\nu^n \cvw1 \nu$, see Lemma \ref{equilibre} in the Appendix below.\\
$\bullet$ The map $\pi : \mathcal{Z} \rightarrow \mathcal{X}$ defined by $\pi(\nu) = \nu_1$ is continuous 
and surjective. \\
$\bullet$ For every $n \geq 1$ and every $\nu = \frac{1}{n} \sum_{i=1}^n \delta_{w_i} \in \mathcal{X}_n$ let $P^n_{\nu}$
be the distribution of 
$T_n = \frac{1}{n} \sum_{i=1}^n \delta_{(w_{\sigma_n(i)},x^n_i)}$
under $\mathbb{P}$. The family
$\Pi = \left\{ P^n_{\nu}, \nu \in \mathcal{X}_n, n \geq 1 \right\}$
of finite measures on the Borel $\sigma$-field on $\mathcal{Z}$ is such that for every $n \geq 1$ and 
every $\nu \in \mathcal{X}_n$ we have $P^n_{\nu}(\pi^{-1}(\{\nu\}^c)=0$.\\
$\bullet$ Let $Q^n$ be the distribution of $\frac{1}{n} \sum_{i=1}^n \delta_{W_i^n}$. For every $n\geq 1$ and every
and every measurable $A \subset \mathcal{M}^{1}_{1}(\mathbb{R}_+ \times \Sigma)$ 
$$\mathbb{P}(\mathcal{V}^n \in A) = \int_{\mathcal{X}_n} P^n_{\nu} (A) Q^n (d \nu).$$

\noindent
All the requirements of Definition 2.1 in \cite{g96} are satisfied by our model thanks to Lemma \ref{lemma2bis}. 
It follows from Theorem 2.3 in \cite{g96} that the sequence $(\mathcal{V}^n)_{n \geq 1}$ obeys an LDP
on $\mathcal{M}^{1}_{1}(\mathbb{R}_+ \times \Sigma)$ with distance $\Delta$ with good rate function

$$
\mathcal{J}(\rho) =  \left\{
\begin{array}{cl}
H(\rho| \rho_1 \otimes \mu) + I^{W}(\rho_1) & \mbox{if \ } \rho_2 = \mu  \\
+ \infty & \mbox{otherwise}.
\end{array}
\right.
$$


\section{Proof of Lemma \ref{clement}}


In order to prove (\ref{varlin}) it is sufficient to consider $\nu,\mu \in M_1(\Sigma)$ such that 
$\mathcal{K}(\nu;\mu) < \infty$ for otherwise the inequality trivially holds. Then  there necessary exists 
(at least) a $\rho \in \mathcal{M}^1_1(\R_+ \times \Sigma)$
such that $\rho_2=\mu$ and $F(\rho)=\nu$. For every such $\rho$ the latter reads
$$\nu(A) = \int_{\R_+ \times A} w \rho_x(dw) \mu(dx)$$
for every measurable $A \subset \Sigma$, hence $\nu \ll \mu$ and 
\begin{equation} \label{filou}
\frac{d\nu}{d\mu}(x) = \int_{\R_+} w \rho_x(dw)
\end{equation}
\noindent
$\mu \ a.s.$ This shows that $\nu \ll \mu$ is a necessary condition for the existence of
$\rho$ such that $\rho_2=\mu$ and $F(\rho)=\nu$ so, as claimed, the announced inequality holds true as an equality if $\nu \ll \mu$
is not satisfied.  Moreover, we have
\begin{eqnarray*}
H(\rho|\rho_1 \otimes \mu) + H(\rho_1 |\xi) & = & H(\rho |\xi \otimes \mu) \\
& =  & \int_{\Sigma} H(\rho_x | \xi ) \mu(dx)
\end{eqnarray*}
\noindent
since $\rho_2=\mu$. It follows from Kullback's inequality (see e.g. Theorem 2.1 in \cite{k38}) that
$$H(\rho_x | \xi ) \geq \Lambda_{\xi}^{\ast}( \int_{\R_+} w \rho_x(dw) )$$
which combined with (\ref{filou}) brings the announced inequality. Actually, it shows that 
a necessary and sufficient condition on $\xi$ to get for every $\nu,\mu \in M_1(\Sigma)$ an equality in (\ref{varlin}) is that 
for every $\nu,\mu \in M_1(\Sigma)$ there exists a 
$\rho \in \mathcal{M}^1_1(\R_+ \times \Sigma)$
such that $\rho_2=\mu$, $F(\rho)=\nu$ and
$$H(\rho_x | \xi ) = \Lambda_{\xi}^{\ast}( \frac{d\nu}{d\mu}(x) )$$ 
$\mu \ a.s.$. According to Kullback's inequality 
this holds true if and only if there exists $\beta_x \in \R$ such that $\mu \ a.s. \ \rho_x$ can take the form
\begin{equation} \label{choucroute}
 \rho_x(dw) = \frac{1}{Z_x} e^{\beta_x w} \xi(dw) 
\end{equation}
\noindent
and still satisfy  
$$\int_{\R_+} w \rho_x(dw)= \frac{d\nu}{d\mu}(x).$$
Let us recall that under (\ref{douane}) the map  $\Lambda_{\xi}$ is $C^{\infty}$ in $\R$, that for every $\alpha \in \R$ 
\begin{eqnarray*}
 \Lambda_{\xi}^{'}(\alpha) & = &  \frac{\int_{\R_+} w e^{\alpha w} \xi(dw)}{\int_{\R_+} e^{\alpha w} \xi(dw)} \\
 & = & \frac{1}{Z} \int_{\R_+} w e^{\alpha w} \xi(dw)
\end{eqnarray*}
\noindent
and that $\Lambda_{\xi}^{''}(\alpha) > 0$ for every $\alpha \in \R$, see Section 2.2.1 in \cite{dz98}.
So finally it appears that 
$$
\lim_{\alpha \rightarrow - \infty} \Lambda_{\xi}^{'}(\alpha) = 0 \mbox{   and   }
\lim_{\alpha \rightarrow + \infty} \Lambda_{\xi}^{'}(\alpha) = + \infty.
$$
is a necessary and sufficient condition to get for every $\nu,\mu \in M_1(\Sigma)$ and every $x \in \Sigma$ such that $\frac{d\nu}{d\mu}(x) > 0$
a $\rho_x$ of the form (\ref{choucroute}) with $\int_{\R_+} w \rho_x(dw)= \frac{d\nu}{d\mu}(x)$ and 
$H(\rho_x | \xi ) = \Lambda_{\xi}^{\ast}( \frac{d\nu}{d\mu}(x))$. If $\frac{d\nu}{d\mu}(x) = 0$ then one takes
$\rho_x=\delta_0$ and still gets $\int_{\R_+} w \rho_x(dw)= \frac{d\nu}{d\mu}(x)$ and 
$H(\rho_x | \xi ) = \Lambda_{\xi}^{\ast}( \frac{d\nu}{d\mu}(x)) = \Lambda_{\xi}^{\ast}(0)$. The announced result follows.


\section{Proof of Theorem \ref{th3}}


In order to prove Theorem \ref{th3} it is sufficient to establish that the distribution of $V^n$
on $\mathcal{M}^{1}_{1}(\mathbb{R}_+ \times \Sigma)$ is a mixture of LDS.
Again we recover the notations of \cite{g96} when identifying the components of the LDS

\noindent
$\bullet$ $\mathcal{Z} = \mathcal{M}^{1}_{1}(\mathbb{R}_+ \times \Sigma)$ is a Polish space when 
endowed with the distance
$\Delta$.\\
$\bullet$ $\mathcal{X} = M_{1}(\Sigma)$ is a Polish space when endowed with the weak convergence topology.\\
$\bullet$ For every $n \geq 1$ we note
$$\mathcal{X}_n = \left\{ \nu \in M^{1}_{1}(\Sigma) \ : \ \exists \ (x_1,\dots,x_n) \in \Sigma^n, \ 
\nu = \frac{1}{n} \sum_{i=1}^n \delta_{x_i} \right\}$$
and according to Varadarajan's Lemma for every $\nu \in \mathcal{X}$ and 
every $n \geq 1$ there exists a $\nu^n \in \mathcal{X}_n$
such that $\nu^n \cvloi \nu$.\\
$\bullet$ The map $\pi : \mathcal{Z} \rightarrow \mathcal{X}$ defined by $\pi(\nu) = \nu_2$ is 
continuous and surjective. \\
$\bullet$ For every $n \geq 1$ and every $\nu = \frac{1}{n} \sum_{i=1}^n \delta_{x_i} \in \mathcal{X}_n$ let $P^n_{\nu}$
be the distribution of 
$T_n = \frac{1}{n} \sum_{i=1}^n \delta_{W^n_i,x_i}$
under $\mathbb{P}$. The family
$\Pi = \left\{ P^n_{\nu}, \nu \in \mathcal{X}_n, n \geq 1 \right\}$
of finite measures on the Borel $\sigma$-field on $\mathcal{Z}$ is such that for every $n \geq 1$ and 
every $\nu \in \mathcal{X}_n$ we have $P^n_{\nu}(\pi^{-1}(\{\nu\}^c)=0$.\\
$\bullet$ Let $Q^n$ be the distribution of $\frac{1}{n} \sum_{i=1}^n \delta_{X_i^n}$. For every $n\geq 1$ and every
and every measurable $A \subset \mathcal{M}^{1}_{1}(\mathbb{R}_+ \times \Sigma)$ 
$$\mathbb{P}(\mathcal{V}^n \in A) = \int_{\mathcal{X}_n} P^n_{\nu} (A) Q^n (d \nu).$$

\noindent
All the requirements of Definition 2.1 in \cite{g96} are satisfied by our model thanks to Theorem \ref{th1}. 
It follows from Theorem 2.3 in \cite{g96} that the sequence $(V^n)_{n \geq 1}$ obeys an LDP
on $\mathcal{M}^{1}_{1}(\mathbb{R}_+ \times \Sigma)$ with distance $\Delta$ with good rate function

$$
J(\rho) =  H(\rho| \rho_1 \otimes \rho_2) + I^{W}(\rho_1) + I^{X}(\rho_2).
$$


\section{Proof of Corollary \ref{contador}}


In view of (\ref{shrek}) a necessary condition on $((W_{i}^{n})_{1 \leq i \leq n})_{n \geq 1}$ for $K=I^X$
is that for every $\nu,\zeta \in M_1(\Sigma)$ and every $((X_{i}^{n})_{1 \leq i \leq n})_{n \geq 1}$
$$I^X(\nu) - I^X(\zeta) \leq \mathcal{K}(\nu,\zeta).$$
If we consider $X^n_1,\dots,X^n_n$ resulting from sampling without replacement on an urn which composition
$x^n_1,\dots,x^n_n$ satisfies $\frac{1}{n} \sum_{i=1}^n \delta_{x^n_i} \cvloi \zeta$ we know that
$\frac{1}{n} \sum_{i=1}^n \delta_{X^n_i}$ satisfies an LDP with good rate function
$$
I^X(\theta) = \left\{
\begin{array}{cl}
0 & \mbox{ if } \theta = \zeta \\
\infty & \mbox{ otherwise. }
\end{array}
\right.
$$
so the announced condition is necessary. It is also clearly sufficient and the announced result follows.


\section{Sample weights Large Deviations Principles} \label{catalogue}


\noindent
In this section we prove the results given in Section \ref{anigo}.

\subsection{Efron's bootstrap and "$m$ out of $n$" bootstrap}

\subsubsection{Proof of Theorem \ref{lemmaE}} \label{recoupling}

The proof relies on the combination of a coupling construction and a Sanov's result in a strong topology.  

\medskip

\noindent
\underline{Coupling Poisson and Multinomial distributions.}

\medskip

\noindent
Let $m,n \geq 1$ and $Z^n_1,\dots,Z^n_n$ be independent random variables such that every $Z^n_i$ is
$\mathcal{P}(m/n)$-distributed. We shall interpret each $Z^n_i$ as the number of balls randomly thrown 
in an urn labeled $i$.
To $(Z^n_1,\dots,Z^n_n)$ we couple $(M^n_1,\dots,M^n_n)$ in the following way
\\
\noindent
$\bullet$ If $\sum_{i=1}^n Z^n_i = m$ we take $M^n_i = Z^n_i$ for every $1 \leq i \leq n$.\\
$\bullet$ If $\sum_{i=1}^n Z^n_i > m$ we pick uniformly at random $\sum_{i=1}^n Z^n_i - m$ balls in the urns. We define 
$(M^n_1,\dots,M^n_n)$ as the new occupation numbers of the urns.\\
$\bullet$ If $\sum_{i=1}^n Z^n_i < m$ we add $m - \sum_{i=1}^n Z^n_i$ balls into the urns. 
The urns are chosen independently 
for each added ball with probability $1 / n$. Again, we denote by 
$(M^n_1,\dots,M^n_n)$ the new occupation numbers of the urns.\\

\begin{lemma} \label{lemmaB} The vector $(M^n_1,\dots,M^n_n)$ is  \hbox{Mult}$_{n}(m,(1 / n,\dots,1 / n))$-distributed.
\end{lemma}

\noindent
The coupling is optimal in the following sense:
\begin{lemma} \label{lemmaA} For every non-negative real numbers $x_1,\dots,x_n$ such 
that $\sum_{i=1}^n x_i = n$ we have
$$W_{1,|\cdot|}(\frac{1}{n} \sum_{i=1}^n \delta_{\frac nm M^n_i},  \frac{1}{n} \sum_{i=1}^n \delta_{\frac nm Z^n_i}) \leq 
W_{1,|\cdot|}(\frac{1}{n} \sum_{i=1}^n \delta_{x_i},  \frac{1}{n} \sum_{i=1}^n \delta_{\frac nm Z^n_i}).$$
\end{lemma}

\bigskip

\noindent
{\bf Proof of Lemma \ref{lemmaB}} Remark that the conditional law of $(Z_1,\ldots,Z_n)$ 
given $\sum_{i=1}^nZ_i^n=k$ is Mult$_n(k,(1/n,\ldots,1/n))$. Thus, if $k=m$ the result is obvious. 
If $k<m$ or $k>m$, the coupling is done such that we move from the law of $k$ balls randomly putted 
in $n$ urns to $m$ balls into the same urns by picking or putting at random. We   finally get the 
Mult$_n(m,(1/n,\ldots,1/n))$ law.
\hfill $\Box$

\bigskip

\noindent
{\bf Proof of Lemma \ref{lemmaA}} By definition
 
\begin{eqnarray*}
W_{1,|\cdot|}(\frac{1}{n} \sum_{i=1}^n \delta_{ \frac nm M^n_i},  \frac{1}{n} \sum_{i=1}^n \delta_{\frac nm Z^n_i})  & = & 
\sup_{f \in C_b (\mathbb{R}_+) \atop ||f||_{L,|\cdot|} \leq 1}
\left\{ \left| \frac{1}{n} \sum_{i=1}^n f(\frac nm M^n_i) - \frac{1}{n} \sum_{i=1}^n f(\frac nm Z^n_i) \right| \right\} \\
& \leq & \frac{1}{m} \sum_{i=1}^n \left| M^n_i -  Z^n_i \right|.
\end{eqnarray*}

\noindent
Due to the construction of $(M^n_1,\dots,M^n_n)$ we have that either all $M^n_i - Z^n_i \leq 0$ 
or all $M^n_i - Z^n_i \geq 0$, so

\begin{eqnarray*}
\frac{1}{m} \sum_{i=1}^n \left| M^n_i - Z^n_i \right| & = & \left| \frac{1}{m} \sum_{i=1}^n  (M^n_i - Z^n_i) \right| \\
& = & \left| 1 - \frac{1}{m} \sum_{i=1}^n   Z^n_i \right|.
\end{eqnarray*}
\noindent
Moreover, for every non-negative real numbers $x_1,\dots,x_n$ such that $\sum_{i=1}^n x_i = n$ we have

\begin{eqnarray*}
W_{1,|\cdot|}(\frac{1}{n} \sum_{i=1}^n \delta_{x_i},  \frac{1}{n} \sum_{i=1}^n \delta_{\frac nm Z^n_i})  & = & 
\sup_{f \in C_b (\mathbb{R}_+) \atop ||f||_{L,|\cdot|} \leq 1}
\left\{ \left| \frac{1}{n} \sum_{i=1}^n f(x_i) - \frac{1}{n} \sum_{i=1}^n f(\frac nm Z^n_i) \right| \right\} \\
& \geq & \left| \frac{1}{n} \sum_{i=1}^n ( x_i - \frac nm Z^n_i) \right| \\
& = & \left| 1 - \frac{1}{m} \sum_{i=1}^n   Z^n_i \right| ,
\end{eqnarray*}
\noindent
which concludes the proof. \hfill $\Box$

\bigskip

\medskip

\noindent
\underline{A Sanov's result in a strong topology.}

\medskip

\noindent
Let $((Z^n_1,\dots,Z^n_n)_{1\leq i \leq n})_{n \geq 1}$ be a triangular array of random variables such that for every
$n \geq 1$ the $Z^n_1,\dots,Z^n_n$ are independent and identically $\mathcal{P}(\lambda_n)$-distributed and
such that $\lim_{n \rightarrow \infty} \lambda_n = \lambda > 0$.

\begin{lemma} \label{evajoly} The sequence 
$(\mathcal{R}^n = \frac{1}{n} \sum_{i=1}^n \delta_{\frac{1}{\lambda_n}Z^n_i})_{n \geq 1}$ obeys an LDP on 
$\mathcal{W}^1(\mathbb{R}_+)$ endowed with the $W_{1,|\cdot|}$ distance with good rate function $H(\cdot|\mathcal{Q}(\lambda))$.
\end{lemma}

\noindent
{\bf Proof of Lemma \ref{evajoly}} According to Theorem 3 in \cite{bj88} the sequence $(\mathcal{R}^n)_{n \geq 1}$ obeys an
LDP on $M_1(\mathbb{R}_+)$ when the latter is endowed with the weak convergence topology with good rate function
$H(\cdot|\mathcal{Q}(\lambda))$. 
Moreover, for every $n \geq 1$ we have $\mathcal{R}^n \in \mathcal{W}^1(\mathbb{R}_+)$ and, since all the exponential moments of
$\mathcal{Q}(\lambda)$ are finite, for every $\nu \in M_1(\mathbb{R}_+)$ if
$H(\nu|\mathcal{Q}(\lambda) < \infty$ then necessarily $\nu \in \mathcal{W}^1(\mathbb{R}_+)$, see e.g. Lemma 2.1 in \cite{es02}.
Thus  $(\mathcal{R}^n)_{n \geq 1}$ obeys an
LDP on $\mathcal{W}^1(\mathbb{R}_+)$ endowed with the weak convergence topology with good rate function
$H(\cdot|\mathcal{Q}(\lambda))$, see Lemma 4.1.5 in \cite{dz98}. 
In particular it is exponentially tight i.e. for every $L> 0$ there exists a compact
$K_L \subset \mathcal{W}^1(\mathbb{R}_+)$ such that
$$\limsup_{n \rightarrow \infty} \frac{1}{n} \log \P( \mathcal{R}^n \nin K_L) < -L.$$
To strengthen the topology for which the LDP holds
it is sufficient to prove that the sequence $(\mathcal{R}^n)_{n \geq 1}$ is exponentially tight in the stronger 
$W_{1,|\cdot|}$ topology, see
Corollary 4.2.6 in \cite{dz98}. To this end we follow the proof of Theorem 1.1. in \cite{www10}. For every $L > 0$ we introduce
$$
A_L= \left\{ \nu \in \mathcal{W}^1(\mathbb{R}_+) : \int_{\R_+} \Lambda^{\ast}_{\mathcal{Q}(\lambda)} (x) \nu(dx) \leq L \right\}.
$$
We prove that for every $L > 0$ the set $C_L = K_L \cap A_L$ is $W_{1,|\cdot|}$-compact. First notice that $A_L$ is 
closed for the weak convergence topology hence $C_L$ is compact for this topology. Thus, to prove
that $C_L$ is $W_{1,|\cdot|}$-compact it remains to prove that
$C_L$ has uniformly integrable first moments. Let $S(N) = \Lambda^{\ast}_{\mathcal{Q}(\lambda)} (N) / N$. According
to Lemma 2.2.20 in \cite{dz98} we have $\lim_{N \rightarrow \infty} S(N) = \infty$. For every $\nu \in C_L$ we have

\begin{eqnarray*}
\int_{x \geq N} x \nu(dx) & \leq & \frac{1}{S(N)} \int_{x \geq N} \frac{\Lambda^{\ast}_{\mathcal{Q}(\lambda)} (x)}{x} x \nu(dx) \\
& \leq & \frac{1}{S(N)} \int_{\R_+} \Lambda^{\ast}_{\mathcal{Q}(\lambda)} (x) \nu(dx) \\
& \leq & \frac{L}{S(N)}
\end{eqnarray*}
\noindent
hence $C_L$ is $W_{1,|\cdot|}$-compact. Since
\begin{equation} \label{abel}
\limsup_{n \rightarrow \infty} \frac{1}{n} \log \P( \mathcal{R}^n \nin C_L) \leq \max \left( -L , 
\limsup_{n \rightarrow \infty} \frac{1}{n} \log \P( \mathcal{R}^n \nin A_L) \right)
\end{equation}
we are led to consider $\P( \mathcal{R}^n \nin A_L)$.
According to Chebichev's inequality for every $\delta \in (0,1)$

\begin{eqnarray}
\P(  \mathcal{R}^n \nin A_L) & = & \P( \frac{1}{n} \sum_{i=1}^n \Lambda^{\ast}_{\mathcal{Q}(\lambda)} (Z^n_i) > L ) \nonumber \\
& \leq 	& \exp(-\delta n L) \left( \E \exp( \delta \Lambda^{\ast}_{\mathcal{Q}(\lambda)} (Z^n_1)) \right)^n. \label{cain} 
\end{eqnarray}

\noindent
But for every $x \in \R$ and every $\lambda > 0$ we have $\Lambda^{\ast}_{\mathcal{Q}(\lambda)}(x)=
\Lambda^{\ast}_{\mathcal{P}(\lambda)}(\lambda x) $ hence
\begin{equation} \label{georges}
\Lambda^{\ast}_{\mathcal{Q}(\lambda)}(x)= \left\{
\begin{array}{cl}
\lambda - \lambda x + \lambda x \log x & \mbox{if } x \geq 0 \\
\infty & \mbox{otherwise.}
\end{array}
\right.
\end{equation}
\noindent
So we get for every $x \in \R$ the relationship
$\Lambda^{\ast}_{\mathcal{Q}(\lambda)}(x) = (\lambda / \lambda_n) \Lambda^{\ast}_{\mathcal{Q}(\lambda_n)}(x)$
so
$$
\P(  \mathcal{R}^n \nin A_L) \leq \exp(-\delta n L) \left( \E \exp( \delta (\lambda / \lambda_n) 
\Lambda^{\ast}_{\mathcal{Q}(\lambda)} (Z^n_1)) \right)^n.$$
According to Lemma 5.1.14 in \cite{dz98} for $n$ large enough we have
$$\E \exp( \delta (\lambda / \lambda_n) 
\Lambda^{\ast}_{\mathcal{Q}(\lambda)} (Z^n_1)) \leq \frac{2}{1-\delta (\lambda / \lambda_n)}$$
which combined with (\ref{abel}) and (\ref{cain}) shows that $(\mathcal{R}^n)_{n \geq 1}$ is exponentially tight in the  
$W_{1,|\cdot|}$ topology. The announced result follows.
 \hfill $\Box$

\medskip

\noindent
\underline{Conclusion of the proof of Lemma \ref{lemmaE}}

\medskip

\noindent
\emph{LD upper bound} Let $F$ be a $W_{1,|\cdot|}$-closed subset of $M^1_1(\R_+)$. Since $M^1_1(\R_+)$
is a closed subset of $\mathcal{W}^1(\R_+)$, $F$ is also a closed subset of $\mathcal{W}^1(\R_+)$ . 
For every $n \geq 1$ let
$Z^n_1,\dots,Z^n_n$ be independent and identically $\mathcal{P}(\lambda_n)$ distributed random variables and
$(M^n_1,\dots,M^n_n)$ the vector associated to $Z^n_1,\dots,Z^n_n$ by the coupling 
construction described in Section \ref{recoupling}.
According to Lemma \ref{lemmaB} $(M^n_1,\dots,M^n_n)$ is Mult$_{n}(m(n),(1 / n,\dots,1 / n))$-distributed, 
which is also the 
distribution of $Z^n_1,\dots,Z^n_n$ conditioned on $\{\sum_{i=1}^n Z^n_i = m(n) \}$. According to 
Lemma \ref{evajoly}

\begin{eqnarray*}
\limsup_{n \rightarrow \infty} \frac{1}{n} \log \P ( \frac{1}{n} \sum_{i=1}^n \delta_{W^n_i} \in F ) & = &
\limsup_{n \rightarrow \infty} \frac{1}{n} \log \P ( \frac{1}{n} \sum_{i=1}^n \delta_{\frac{n}{m(n)} M^n_i} \in F ) \\
 & \leq & 
\limsup_{n \rightarrow \infty} \frac{1}{n} \log \P ( \frac{1}{n} \sum_{i=1}^n \delta_{\frac{1}{\lambda_n} Z^n_i} \in F ) \\
& & \ \ \ \ - \liminf_{n \rightarrow \infty} \frac{1}{n} \log \P (\sum_{i=1}^n Z^n_i = m(n)) \\
& \leq & - \inf_{\nu \in F } H(\nu | \mathcal{Q}(\lambda))
\end{eqnarray*}

\noindent 
since $\liminf_{n \rightarrow \infty} \frac{1}{n} \log \P (\sum_{i=1}^n Z^n_i = m(n)) = 0$.

\bigskip

\noindent
\emph{LD lower bound} Let $O$ be an open subset of $M^1_1(\R_+)$. In order to prove the 
LD lower bound it is sufficient to prove
that for every $\rho \in  O$ we have

$$\liminf_{n \rightarrow \infty} \frac{1}{n} \log \P ( \frac{1}{n} \sum_{i=1}^n \delta_{W^n_i} \in O ) \geq - 
H(\rho | \mathcal{Q}(\lambda)).$$
So let $\rho \in O$ and $\varepsilon > 0$ be small enough to ensure that $B(\rho,\varepsilon) \subset O$. 
We can assume that
$H(\nu | \mathcal{Q}(\lambda)) < \infty$ for otherwise the LD lower bound trivially holds. According 
to Lemma \ref{equilibre} proved in the Appendix for every $\rho \in M^1_1(\R_+)$ there exists 
a triangular array of non-negative real numbers 
$((x^n_i)_{1 \leq i \leq n})_{n \geq 1}$ such that $\rho^n = \frac{1}{n} \sum_{i=1}^n \delta_{x^n_i} \cvw1 \rho$ 
and such that for every $n \geq 1$ we have $\sum_{i=1}^n x^n_i = n$.
\\

\noindent
Again, let $Z^n_1,\dots,Z^n_n$ be independent and identically $\mathcal{P}(\lambda_n)$ distributed random variables and
$(M^n_1,\dots,M^n_n)$ the vector associated to $Z^n_1,\dots,Z^n_n$ by the coupling construction described above.
In particular

$$W_{1,|\cdot|}(\frac{1}{n} \sum_{i=1}^n \delta_{\frac{n}{m(n)} M^n_i},\frac{1}{n} \sum_{i=1}^n \delta_{\frac{1}{\lambda_n} Z^n_i}) 
\leq W_{1,|\cdot|}(\frac{1}{n} \sum_{i=1}^n \delta_{\frac{1}{\lambda_n} Z^n_i}, \frac{1}{n} \sum_{i=1}^n \delta_{x^n_i}).$$
So,
\begin{eqnarray*}
\mathbb{P}(\frac{1}{n} \sum_{i=1}^n \delta_{ W^n_i} \in O) & = & 
\mathbb{P}(\frac{1}{n} \sum_{i=1}^n \delta_{\frac{n}{m(n)} M^n_i} \in O)	\\
& \geq &  \mathbb{P}( W_{1,|\cdot|} (\frac{1}{n} \sum_{i=1}^n \delta_{\frac{n}{m(n)} M^n_i}, 
\frac{1}{n} \sum_{i=1}^n \delta_{\frac{1}{\lambda_n} Z^n_i}) 
< \varepsilon / 3, \\
			& & 						W_{1,|\cdot|}(\rho^n,\frac{1}{n} \sum_{i=1}^n \delta_{\frac{1}{\lambda_n} Z^n_i}) < \varepsilon / 3, \ 
													  W_{1,|\cdot|}(\rho^n,\rho) < \varepsilon / 3 ) \\
				 & \geq & \mathbb{P}(  W_{1,|\cdot|}(\rho^n,\frac{1}{n} \sum_{i=1}^n \delta_{\frac{1}{\lambda_n} Z^n_i}) < \varepsilon / 3, \ 
													  W_{1,|\cdot|}(\rho^n,\rho) < \varepsilon / 3 ) \\
				& \geq & \mathbb{P}(  W_{1,|\cdot|}(\rho,\frac{1}{n} \sum_{i=1}^n \delta_{\frac{1}{\lambda_n} Z^n_i}) < \varepsilon / 6 )
\end{eqnarray*}
\noindent
for every $n$ large enough and the LD lower bound follows from Lemma \ref{evajoly}. \hfill $\Box$

\subsubsection{Proof of Corollary \ref{nkoulou}}

First notice that $\Lambda_{\mathcal{Q}(\lambda)}(\alpha) = \lambda (e^{\alpha / \lambda} - 1)$ so condition (\ref{gargamel}) 
is satisfied. Using (\ref{georges}) we immediately obtain the announced result.  \hfill $\Box$

\hfill $\Box$

\subsubsection{Proof of Corollary \ref{amalfitano}}
Let us prove the inequality in Corollary \ref{amalfitano}.
Let $\nu \in M_1(\Sigma)$ be such that $\inf_{\eta \in \mathcal{Z}} H(\nu|\eta) < \infty$. Necessarily for every $\eta \in \mathcal{Z}$
such that $H(\nu|\eta) < \infty$ we have $\nu \ll \eta$. Now
for every such $\eta$ consider 
$\rho \in \mathcal{M}_1^1 (\mathbb{R}_+ \times \Sigma)$ defined by $\rho_2 = \eta$ while the regular conditional distribution 
of its first marginal given the second is
$$\rho_x(dw)=
\left\{
\begin{array}{cl}
\mathcal{F}(\lambda,\lambda \frac{d\nu}{d\eta}(x))(dw) & \mbox{ if  } \frac{d\nu}{d\eta}(x) > 0 \\
\delta_0(dw) & \mbox{ otherwise.}
\end{array}
\right.
$$
Let us check that $F(\rho)=\nu$. Indeed, every measurable $A \subset \Sigma$  can be decomposed
into $A = A_{\nu} \cup A_{\nu}^{\bot}$ where $ A_{\nu} = A \cap \mbox{Support}(\nu)$ and
\begin{eqnarray*}
F(\rho)(A) & = & \int_{\R_+ \times A} w \rho_x(dw) \eta(dx) \\
& = & \int_{A_{\nu}} \left( \int_{\R_+} w \mathcal{F}(\lambda,\lambda \frac{d\nu}{d\eta}(x))(dw) \right) \eta(dx) + 
 \int_{A_{\nu}^{\bot}} \left( \int_{\R_+} w \delta_0(dw) \right) \eta(dx) \\
& = & \int_{A_{\nu}} \frac{d\nu}{d\eta}(x) \eta(dx) \\
& = & \nu(A_{\nu}) = \nu(A).
\end{eqnarray*}

\noindent
Moreover by taking $\Sigma_1=\Sigma \cap \{\frac{d\nu}{d\eta}(x) > 0 \}$ and
 $\Sigma_2=\Sigma \cap \{\frac{d\nu}{d\eta}(x) = 0 \}$ we get
\begin{eqnarray*}
H(\rho|\rho_1 \otimes \eta) + H(\rho_1 | \mathcal{Q}(\lambda)) + I^X(\eta) & = & H(\rho | \mathcal{Q}(\lambda) \otimes \eta) \\
& = & \int_{\Sigma} H( \rho_x(\cdot)|\mathcal{Q}(\lambda)) \eta(dx) \\
& = & \int_{\Sigma_1} H( \mathcal{P}(\lambda \frac{d\nu}{d\eta}(x))|\mathcal{P}(\lambda)) \eta(dx) +\\
&   & \ \ \ \ \ \ + \int_{\Sigma_2} H( \delta_0|\mathcal{P}(\lambda)) \eta(dx) \\
& = & \int_{\Sigma_1} \lambda (1 - \frac{d\nu}{d\eta}(x) + \frac{d\nu}{d\eta}(x) \log \frac{d\nu}{d\eta}(x) ) \eta(dx) + \lambda \eta(\Sigma_2) \\
& = & \lambda H(\nu|\eta),
\end{eqnarray*}
\noindent
hence $K(\nu) \leq \lambda \inf_{\eta \in \mathcal{Z}} H(\nu|\eta)$. \hfill $\Box$

\subsection{Iid weighted bootstrap} The following

\begin{lemma} \label{gcontinuite}
$\mathcal{G}$ is continuous.
\end{lemma}

\noindent 
is the main argument in the proof of Theorem \ref{pgdiid}

\subsubsection{Proof of Lemma \ref{gcontinuite}}
We consider a sequence of probability measures
$(\rho^n)_{n \geq 1}$ such that $\rho^n \cvw1 \rho$ and a sequence of positive numbers $m_n \rightarrow m > 0$ and prove
that $\mathcal{G}(\rho^n,m_n) \cvw1 \mathcal{G}(\rho,m)$. It is not difficult to show that for every 
$f \in C_b(\mathbb{R}_{+}^{*})$, every $\rho \in \mathcal{W}^1(\mathbb{R}_+)$ and every $m>0$ we have

$$\int_{\mathbb{R}_{+}} f(x) \mathcal{G}(\rho,m)(dx) = \int_{\mathbb{R}_{+}} f(\frac{x}{m}) \rho(dx).$$ 
So $\int_{\mathbb{R}_{+}} x \mathcal{G}(\rho^n,m_n)(dx)  =  \int_{\mathbb{R}_{+}} \frac{x}{m_n} \rho^n (dx)
			 \rightarrow  \frac{1}{m} \int_{\mathbb{R}_{+}} x \rho (dx) 
			 = \int_{\mathbb{R}_{+}} x \mathcal{G}(\rho,m)(dx)$. Hence, according to Theorem 7.11 
in \cite{v03} we are left to prove that for  every $f \in C_b(\mathbb{R}_{+}^{*})$ we have
$$ \int_{\mathbb{R}_{+}} f(x) \mathcal{G}(\rho^n,m_n)(dx)  \rightarrow 
\int_{\mathbb{R}_{+}} f(x) \mathcal{G}(\rho,m)(dx).$$
Since $\rho^n \cvw1 \rho$ we have $\rho^n \cvloi \rho$ hence for every $\varepsilon > 0$ there exists a compact 
$K_{\varepsilon} \subset \mathbb{R}_{+}^{*}$ such that 
$\sup_{n \geq 1} \{ \rho^n (K_{\varepsilon}^c), \rho (K_{\varepsilon}^c) \} < \varepsilon$. Thus
\begin{eqnarray}
& & \left| \int_{\mathbb{R}_{+}} f(x) \mathcal{G}(\rho^n,m_n)(dx) - \int_{\mathbb{R}_{+}} 
f(x) \mathcal{G}(\rho,m)(dx) \right|  \nonumber \\
& = &
\left| \int_{\mathbb{R}_{+}} f(\frac{x}{m_n}) \rho^n(dx) - \int_{\mathbb{R}_{+}} f(\frac{x}{m}) \rho(dx) \right| \nonumber \\
& \leq & \left| \int_{K_{\varepsilon}} f(\frac{x}{m_n}) \rho^n(dx) - 
\int_{K_{\varepsilon}} f(\frac{x}{m}) \rho(dx) \right| + 2 \varepsilon ||f||_{\infty} \nonumber \\
& \leq & \int_{K_{\varepsilon}} \left| f(\frac{x}{m_n}) - f(\frac{x}{m}) \right| \rho^n(dx)  + \nonumber \\
& & \ \ \ \ \ \ \ 
 + \left| \int_{K_{\varepsilon}} f(\frac{x}{m}) \rho^n(dx) - \int_{K} f(\frac{x}{m}) \rho(dx) \right| 
 + 2 \varepsilon ||f||_{\infty}  \label{hubertnapasletemps}
\end{eqnarray}
Since $f$ is continuous it is uniformly continuous over compact 
subsets of $\mathbb{R}_{+}^{*}$ so there is an $N_0$ such that 
for every $n \geq N_0$ 
$$\sup_{x \in K_{\varepsilon}} \left| f(\frac{x}{m_n}) - f(\frac{x}{m}) \right| < \varepsilon.$$
By taking $N_0$ large enough we can make the second term in (\ref{hubertnapasletemps})
 as small as desired, which concludes the proof.
\hfill $\Box$

\subsubsection{Proof of Theorem \ref{pgdiid}}
According to \cite{www10}, that all the exponential 
moments of $\xi$ are finite is a 
necessary and sufficient condition
for $\frac{1}{n} \sum_{i=1}^n \delta_{Y_i}$ to satisfy an LDP on $\mathcal{W}^1(\mathbb{R}_+)$ endowed with $W_{1,|\cdot|}$ with good 
rate function $H(\ \cdot \ |\mu)$. Next, since
$$
\begin{array}{rccc}
G: & \mathcal{W}^{1}(\mathbb{R}_+) & \rightarrow & \mathbb{R}_+ \\
   &  \nu & \mapsto & \int_{\mathbb{R}_+} x \nu(dx)
\end{array}
$$
is continuous when $\mathcal{W}^{1}(\mathbb{R}_+)$ is furnished with the $W_1$ 
distance we obtain that $(\frac{1}{n} \sum_{i=1}^n \delta_{Y_i},
\frac{1}{n} \sum_{i=1}^n Y_i)$ satisfies an LDP on the product 
space $\mathcal{W}^{1}(\mathbb{R}_+) \times \mathbb{R}_+$ with good rate function

$$
R(\rho,x) =  \left\{
\begin{array}{cl}
H(\rho| \mu_1) & \mbox{if \ } \int_{\mathbb{R}_+} u \rho(du) = x  \\
+ \infty & \mbox{otherwise}.
\end{array}
\right.
$$
Since $\Lambda^*(0)=\infty$ we have $R(\rho,0)= \infty$. Hence, according to Lemma 4.1.5 in \cite{dz98}, an LDP 
 for $(\frac{1}{n} \sum_{i=1}^n \delta_{Y_i},
\frac{1}{n} \sum_{i=1}^n Y_i)_{n \geq 1}$ holds on $\mathcal{W}^{1}(\mathbb{R}_+) \times \mathbb{R}_{+}^{*}$ with 
the same rate function $R$. Since $\mathcal{G}$ is continuous and
$$\frac{1}{n} \sum_{i=1}^n \delta_{W_i^n} = \mathcal{G}(\frac{1}{n} 
\sum_{i=1}^n \delta_{Y_i}, \frac{1}{n} \sum_{i=1}^n Y_i)$$
$(\mathcal{S}^n)_{n \geq 1}$ obeys an LDP on $\mathcal{W}^{1}(\mathbb{R}_+)$ due
to the contraction principle, see Theorem 4.2.1 in \cite{dz98}. Finally, the announced result follows from
Lemma 4.1.5 in \cite{dz98} since for every $n \geq 1 \  \mathcal{S}^n \in M_1^1(\mathbb{R}_+)$ and the latter is a closed subset of
$\mathcal{W}^{1}(\mathbb{R}_+)$.
\hfill $\Box$

\subsubsection{Proof of Corollary \ref{vockler}}
Let $\nu \in M_1(\Sigma)$ be such that $\nu \ll \mu$ for otherwise we already know that $\mathcal{K}(\nu;\mu) = + \infty$. 
We have

\begin{eqnarray*}
\mathcal{K}(\nu;\mu) & = & \inf_{\rho_{x} : F(\rho_{x}\otimes \mu)=\nu } \left\{
H(\rho |\rho_1 \otimes \mu) + \inf_{m > 0} \left\{ H( \mathcal{G}(\frac{1}{m},\rho_1 )| \xi ) \right\} \right\} \\
& = & \inf_{m > 0} \inf_{\rho_{x} : F(\rho_{x}\otimes \mu)=\nu } \left\{
H(\rho |\rho_1 \otimes \mu) + H(\rho_1 | \mathcal{G}(m,\xi ) \right\} \\
& = & \inf_{m > 0} \inf_{\rho_{x} : F(\rho_{x}\otimes \mu)=\nu } \left\{ H(\rho |\mathcal{G}(m,\xi ) \otimes \mu) \right\} \\
& = & \inf_{m > 0} \inf_{\rho_{x} : F(\rho_{x}\otimes \mu)=\nu } \left\{ \int_{\Sigma} H(\rho_x |\mathcal{G}(m,\xi )) \mu(dx) \right\} \\
& \geq & \inf_{m > 0} \int_{\Sigma} \Lambda^{\ast}_{\xi}( m \frac{d\nu}{d\mu}(x)) \mu(dx)
\end{eqnarray*}

\noindent
where, to establish the last inequality and characterize the equality, we proceed as in the proof of Lemma \ref{clement}.
\hfill $\Box$

\subsection{The multivariate hypergeometric bootstrap}
Since the $W^n_i$'s take values in the finite set $\{0,\dots,K\}$ the proof of this LDP closely 
follows the proof of Sanov's Theorem for finite alphabets as exposed in
Section 2.1.1 in \cite{dz98} . This is the reason why 
we shall not give all the details of the proof. The main step is to notice that 
for every $(w^{n}_{1},\dots,w^{n}_{n}) \in \{0,\dots,K\}^n$ such that $\sum_{i=1}^n w^n_i = n$  if we introduce
$$(\nu^n(1),\dots,\nu^n(K)) = (\frac{1}{n} \sum_{i=1}^n \delta_{w^{n}_{i}}(1),\dots,\frac{1}{n} 
\sum_{i=1}^n \delta_{w^{n}_{i}}(K))$$
we get $\sum_{k=1}^K k \nu^n(k) =1$ and
\begin{equation}
\mathbb{P}(W^{n}_{1}=w^{n}_{1},\dots,W^{n}_{n}=w^{n}_{n}) = R(n)
\Pi_{k=1}^K \left( C^k_K (\frac{1}{K})^k (1 - \frac{1}{K})^{K-k} \right)^{n \nu^n(k)} \label{stirling}
\end{equation}
\noindent
with $\frac{1}{n} \log R(n) \rightarrow 0$. For every integer $n \geq 1$, every $\nu \in M^1_1(\{0,\dots,K\})$ 
of the form $(\frac{k_1}{n},\dots,\frac{k_K}{n})$
with $k_1,\dots,k_K$ integers there exists $T_n(\nu)$ 
vectors $(w^{n}_{1},\dots,w^{n}_{n}) \in \{0,\dots,K\}^n$ such that
$\nu = \frac{1}{n} \sum_{i=1}^n \delta_{w^{n}_{i}}$ with 
\begin{equation}
(n+1)^{-K} e^{-nH(\nu)} \leq T_n(\nu) \leq  e^{-nH(\nu)} \label{alphabet}
\end{equation}
\noindent
where $H(\nu)= - \sum_{i=1}^K \nu(i) \log \nu(i)$, see Lemma 2.1.8 in \cite{dz98}. Hence, combining (\ref{stirling}) and
(\ref{alphabet}) for every integer 
$n \geq 1$, every $\nu \in M^1_1(\{0,\dots,K\})$ of the form $(\frac{k_1}{n},\dots,\frac{k_K}{n})$
with $k_1,\dots,k_K$ integers we get

$$R_1(n) e^{-nH(\nu|\mathfrak{B}(K,\frac{1}{K}))} \leq \mathbb{P}( \mathcal{S}^n = \nu ) \leq 
R_2(n) e^{-nH(\nu|\mathfrak{B}(K,\frac{1}{K}))}$$
with $\frac{1}{n} \log R_{1,2}(n) \rightarrow 0$. The proof then follows as in Theorem 2.1.10 in \cite{dz98} 
until its conclusion.
\hfill $\Box$

\subsection{A bootstrap generated from deterministic weights}

\subsubsection{Proof of Corollary \ref{cavendish}} 
First we consider $\alpha > 0$.
Let $\nu \in M_1(\Sigma)$ be such that $\frac{\mu - (1 -  \alpha) \nu}{\alpha} \in M_1(\Sigma)$. Then $\rho  
\in \mathcal{M}^1_1(\R_+ \times \Sigma)$ 
defined by $\rho_1=\gamma$ while the regular conditional distribution of its second marginal given the first is
$$
\rho_{\frac{1}{1-\alpha}}(dx) = \nu(dx) \mbox{   and   }  \rho_0(dx) = \frac{\mu(dx) - (1 -  \alpha) \nu(dx)}{\alpha}
$$
is the only element of $\mathcal{M}^1_1(\R_+ \times \Sigma)$ that
satisfies $F(\rho) = \nu, \rho_1 = \gamma$ and $\rho_2=\mu$. Since for this particular $\rho$ we have
$$
H(\rho | \rho_1 \otimes \rho_2) = (1-\alpha)H(\nu|\mu) + \alpha H \left( \frac{\mu - (1 -  \alpha) \nu}{\alpha} |\mu \right )
$$
we obtain an upper-bound on $\mathcal{K}$ as announced. To prove the reverse inequality let $\nu \in M_1(\Sigma)$ be such that
$\mathcal{K}(\nu;\mu) < \infty$ for otherwise the announced result trivially holds. Necessarily there exists
a $\rho  
\in \mathcal{M}^1_1(\R_+ \times \Sigma)$ such that $F(\rho) = \nu, \rho_1 = \gamma$ and $\rho_2=\mu$. These conditions are only met
by the probability measure $\rho$ introduced above. In particular $\frac{\mu - (1 -  \alpha) \nu}{\alpha}$ must
be a probability and the reverse inequality holds.

\bigskip

\noindent 
If $\alpha=0$ then $\gamma=\delta_1$ so the only $\nu \in M_1(\Sigma)$ such that there exists
$\rho \in \mathcal{M}^1_1(\R_+ \times \Sigma)$ such that $F(\rho) = \nu$ and $\rho_2=\mu$ is $\mu$ and necessarily
$\rho= \delta_1 \otimes \mu$. The announced result follows.
\hfill $\Box$

\subsubsection{Proof of Corollary \ref{hinault}} First we consider $\alpha > 0$.
Let $\nu \in M_1(\Sigma)$ be such that $\inf_{\zeta \in  \mathcal{E}_{\nu}} \mathcal{U}(\nu,\zeta) < \infty$ holds. Then
for every $\zeta \in \mathcal{E}_{\nu}$ we have $\frac{\zeta - (1 -  \alpha) \nu}{\alpha} \in M_1(\Sigma)$ and 
$\rho \in \mathcal{M}^1_1(\R_+ \times \Sigma)$ defined by $\rho_1=\gamma$ and 
$$
\rho_{\frac{1}{1-\alpha}}(dx) = \nu(dx) \mbox{   and   }  \rho_0(dx) = \frac{\zeta(dx) - (1 -  \alpha) \nu(dx)}{\alpha}
$$
satisfies $F(\rho) = \nu, \rho_1 = \gamma$ and $\rho_2=\zeta$ hence $\mathcal{K}(\nu;\mu) \leq \mathcal{U}(\nu,\zeta)$
and the announced upper bound on $\mathcal{K}$ follows.To prove the reverse inequality let us assume that 
$\mathcal{K}(\nu;\mu) < \infty$. Then there exists a $\rho \in \mathcal{M}^1_1(\R_+ \times \Sigma)$ such that
$F(\rho) = \nu$ and $\rho_1 = \gamma$. Necessarily $\rho_{\frac{1}{1-\alpha}}(dx) = \nu(dx)$ and $\rho_2$ is such that
$\rho_0(dx) = \frac{\zeta(dx) - (1 -  \alpha) \nu(dx)}{\alpha} \in M_1(\Sigma)$. Thus $\mathcal{E}_{\nu}$ is non-empty
and

\begin{eqnarray*}
&&H(\rho | \rho_1 \otimes \rho_2) + I^W(\rho_1) +  I^X(\rho_2)\\
 & = & (1-\alpha)H(\nu|\rho_2) + \alpha H \left( \frac{\rho_2 - (1 -  \alpha) \nu}{\alpha} |\rho_2 \right ) I^X(\rho_2)\\
 & \geq & \inf_{\rho_2 \in  \mathcal{E}_{\nu}} \mathcal{U}(\nu,\rho_2)
\end{eqnarray*}

\bigskip

\noindent 
If $\alpha=0$ then $\gamma=\delta_1$ and for every $\nu \in M_1(\Sigma)$ there is only one 
$\rho \in \mathcal{M}^1_1(\R_+ \times \Sigma)$ such that $F(\rho) = \nu$ which is $\rho=\delta_1 \otimes \nu$ 
and the announced result immediately follows.
\hfill $\Box$
\subsection{The $k$-blocks bootstrap}
\subsubsection{Proof of Theorem \ref{th:kk}}
Let us apply Corollary \ref{nkoulou} to $\mathcal L^n=\frac 1n \sum_{i=1}^{n }W_i^n\delta_{(x_{i}^n,\ldots,x_{i+k-1}^n)}$. 
We obtain the desired result using the contraction principle on the map
$$
\begin{array}{rccc}
\mathcal{H}: & M_1(\Sigma^k) & \rightarrow &M_1(\Sigma )  \\
   &  \rho^{(k)}& \mapsto & \mathcal{H}(\rho^{(k)})=\frac1k \sum_{i=1}^k  \rho^{(k)}_i
\end{array}
$$
which is continuous because it is Lipschitz.\hfill$\Box$

\subsubsection{Proof of Corollary \ref{nasri}}

Let $\nu \in M_1(\Sigma)$. Since $\mathcal{H}(\nu^{\otimes k}) = \nu$ and in this particular case 
$\mu^{(k)}= \mu^{\otimes k}$ we get 

\begin{eqnarray*}
\widetilde{\mathcal K}(\nu;\mu^{(k)}) & \leq & \frac{1}{k} H(\nu^{\otimes k} | \mu^{\otimes k}) \\
 & = & H(\nu | \mu).
\end{eqnarray*}
\noindent
On the other hand for every $\rho^{(k)}$ such that $\mathcal{H}(\rho^{(k)})=\nu$ we have

$$
 H(\nu | \mu)  =  H(\frac1k \sum_{i=1}^k  \rho^{(k)}_i | \mu) 
  \leq  \frac1k \sum_{i=1}^k H(\rho^{(k)}_i | \mu) 
$$
\noindent
since $H(\cdot|\mu)$ is convex. To conclude the proof just notice that
\begin{eqnarray*}
 H(\rho^{(k)} | \mu^{\otimes k}) & = & H(\rho^{(k)} | \otimes_{i=1}^k \rho^{(k)}_i) 
+ H( \otimes_{i=1}^k \rho^{(k)}_i | \mu^{\otimes k}) \\
& \geq & H( \otimes_{i=1}^k \rho^{(k)}_i | \mu^{\otimes k}) \\
& = & \sum_{i=1}^k H(\rho^{(k)}_i | \mu). 
\end{eqnarray*}
\noindent
\hfill$\Box$


\appendix

\section{Topological properties of $(\mathcal{M}^{1}_{1}(\mathbb{R}_{+} \times \Sigma),\Delta)$ and 
$(M^{1}_{1}(\mathbb{R}_+),W_1)$}


\subsection{$(\mathcal{M}^{1}_{1}(\mathbb{R}_{+} \times \Sigma),\Delta)$ is a Polish space}

\subsubsection{$(\mathcal{M}^{1}_{1}(\mathbb{R}_{+} \times \Sigma),\Delta)$ is complete}
Let $(\rho^n)_{n \geq 1}$ be a Cauchy sequence of 
elements of $(\mathcal{M}^{1}_{1}(\mathbb{R}_+ \times \Sigma),\Delta)$.
In particular $(\rho^n_1)_{n \geq 1}$ is a Cauchy sequence of elements of $(M^{1}_{1}(\mathbb{R}_+),W_1)$
which is a complete space (see e.g. \cite{bo08})
so there exists a $\rho_1 \in M^{1}_{1}(\mathbb{R}_{+})$ such that $W_1(\rho^n_1,\rho_1) \rightarrow 0$. 
Furthermore $(\rho^n)_{n \geq 1}$ is a Cauchy sequence of elements 
of $(M_{1}(\mathbb{R}_+ \times \Sigma),\beta_{BL,\delta})$
hence there exists a $\gamma \in  M_{1}(\mathbb{R}_{+} \times \Sigma)$ such 
that $\beta_{BL,\delta}(\rho^n,\gamma) \rightarrow 0$. Necessarily
we have $\gamma_1 = \rho_1$ hence $\Delta (\rho^n,\gamma) \rightarrow 0$.

\subsubsection{$(\mathcal{M}^{1}_{1}(\mathbb{R}_{+} \times \Sigma),\Delta)$ is separable}
Let $R$ and $E$ be dense countable subsets of $\mathbb{R}_+$ and $\Sigma$ respectively. It is sufficient to prove that
$$\bigcup_{n \geq 1} \left\{ \frac{1}{n} \sum_{i=1}^n \delta_{(u^n_i,v^n_i)} \ : \ 
\mbox{for every } 1 \leq i \leq n, \ u^n_i \in R \ \mbox{and }v^n_i \in E \right\}$$
is dense in $\mathcal{N}_1^1(\mathbb{R}_{+} \times \Sigma)$
endowed with the distance $\Delta$. Indeed, since $\mathcal{M}_1^1(\mathbb{R}_{+} \times \Sigma)$
is a closed subset of $\mathcal{N}_1^1(\mathbb{R}_{+} \times \Sigma)$ the announced claim follows.
So let $\rho \in \mathcal{N}_1^1(\mathbb{R}_{+} \times \Sigma)$. In particular
$\int_{\R_+} x \rho_1 (dx) < \infty$. Lets us denote by
$(X_1,Y_1),\dots,(X_n,Y_n),\dots$ a sequence of independent random variables with common
distribution $\rho$. According to Varadarajan's Lemma

\begin{equation}
\frac{1}{n} \sum_{i=1}^n \delta_{(X^n_i,Y^n_i)} \cvloi \rho \mbox{\ \ almost surely}
\label{lutin}
\end{equation}
and according to the strong law of Large Numbers
\begin{equation}
\frac{1}{n} \sum_{i=1}^n X^n_i \rightarrow \int_{\mathbb{R}_+} x \rho_1(dx) \mbox{\ \ almost surely}
\label{oiseau}
\end{equation}
So there exists a family $((x^n_i,y^n_i)_{1 \leq i \leq n})_{n \geq 1}$ of elements of 
$\mathbb{R}_{+} \times \Sigma$ such that $$\frac{1}{n} \sum_{i=1}^n \delta_{(x^n_i,y^n_i)} \cvloi \rho 
\mbox{\ \ and \ \ \ } \frac{1}{n} \sum_{i=1}^n X^n_i \rightarrow \int_{\mathbb{R}_+} x \rho_1(dx).$$
Since $R \times E$ is dense in $\mathbb{R}_{+} \times \Sigma$ for every $n \geq 1$ 
and every $1 \leq i \leq n$ there exists
$(u^n_i,v^n_i) \in R \times E$ such that $\max \{ |u_i^n - x_i^n|, d_{\Sigma}(v^n_i,y^n_i) \} \leq 2^{-n}.$
Clearly,
$$W_{1,\tilde{d}_{2,+}}(\frac{1}{n} \sum_{i=1}^n \delta_{(u^n_i,v^n_i)},\frac{1}{n} 
\sum_{i=1}^n \delta_{(x^n_i,y^n_i)}) \leq 2^{-n}$$
hence
$W_{1,\tilde{d}_{2,+}}(\frac{1}{n} \sum_{i=1}^n \delta_{(u^n_i,v^n_i)},\rho) \rightarrow 0$.
Since both $\beta_{BL,\delta}$ and $W_{1,\tilde{d}_{2,+}}$ metrize the weak convergence topology we have that 
$$\beta_{BL,\delta}(\frac{1}{n} \sum_{i=1}^n \delta_{(u^n_i,v^n_i)},\rho) \rightarrow 0.$$
Furthermore we have $\frac{1}{n} \sum_{i=1}^n \delta_{u^n_i} \cvloi \rho_1$ and since
\begin{eqnarray*}
\left| \frac{1}{n} \sum_{i=1}^n u^n_i - \frac{1}{n} \sum_{i=1}^n x^n_i \right|  & \leq & 
\frac{1}{n} \sum_{i=1}^n |u^n_i -x^n_i| \\
&  \leq & 2^{-n} 
\end{eqnarray*}
\noindent
we obtain that $\frac{1}{n} \sum_{i=1}^n u^n_i \rightarrow \int_{\mathbb{R}_+} x \rho_1(dx)$ which 
together with the latter weak convergence
ensure that $\frac{1}{n} \sum_{i=1}^n \delta_{u^n_i} \cvw1 \rho_1$, according to Theorem 7.11 in \cite{v03}.

\subsection{An approximation result on $(M^{1}_{1}(\mathbb{R}_+),W_1)$.}

\begin{lemma} \label{equilibre}
For every $\rho \in M^{1}_{1}(\mathbb{R}_+)$ there exists for every $n \geq 1$ an elements $\rho^n$ of
$$
\mathcal{A}_n = \left\{ \nu \in M^{1}_{1}(\mathbb{R}_+), \ 
\exists \ (x_1,\dots,x_n) \in (\mathbb{R}_+)^n, \nu = \frac{1}{n} 
\sum_{i=1}^n \delta_{x_n} \right\}
$$
such that $\nu^n \cvw1 \rho$. 
\end{lemma}
\noindent
\textbf{Proof} By the same kind of argument as in the separability proof above one can construct a sequence
$(\frac{1}{n} \sum_{i=1}^n \delta_{u^n_i})_{n \geq 1}$ such that $\frac{1}{n} \sum_{i=1}^n \delta_{u^n_i} \cvw1  \rho$.
In particular $\frac{1}{n} \sum_{i=1}^n u^n_i \rightarrow 1$. So we only need to modify the $u_i$'s in such a way that 
for every $n \geq 1$ their total sum equals $n$. For a fixed $n$ we have three possibilities:
\noindent
\\
$\bullet$ \underline{ If $\sum_{i=1}^n u_i = n$} we take $w_i^n = u_i$ for every $1 \leq i \leq n$.\\
$\bullet$ \underline{ If $\sum_{i=1}^n u_i > n$} we look at the $u_i$'s as the occupation 
masses of $n$ cells by a mass $u_i$ each.
We pick uniformly at random the excess of mass until we get new occupation masses $w_1^n,\dots,w_n$ such 
that $\sum_{i=1}^n w^n_i = n$.\\
$\bullet$ \underline{ If $\sum_{i=1}^n u_i < n$} again we look at the $u_i$'s as the occupation masses
and add mass uniformly at random into
the $n$ cells until they contain a total mass of $n$. We call $w_1^n,\dots,w_n$ the final occupation masses.

\noindent
In all the cases considered above we have
\begin{eqnarray*}
W_1(\frac{1}{n} \sum_{i=1}^n \delta_{u^n_i},\frac{1}{n} \sum_{i=1}^n \delta_{w_i^n}) & = & 
\sup_{f \in C_b (\mathbb{R}_+) \atop ||f||_{L} \leq 1}
\left\{ \left| \frac{1}{n} \sum_{i=1}^n f(u^n_i) - \frac{1}{n} \sum_{i=1}^n f(w^n_i) \right| \right\} \\
	& \leq & \sup_{f \in C_b (\mathbb{R}_+) \atop ||f||_{L} \leq 1} \frac{1}{n} \sum_{i=1}^n | f(u^n_i) - f(w^n_i) | \\
	& \leq & \frac{1}{n} \sum_{i=1}^n | u^n_i - w^n_i | \\
	& = & \left| \frac{1}{n} \sum_{i=1}^n  u^n_i - \frac{1}{n} \sum_{i=1}^n w^n_i \right| \\
	& = & \left| \frac{1}{n} \sum_{i=1}^n  u^n_i - 1 \right| \rightarrow 0.
\end{eqnarray*}  \hfill $\Box$
\bibliographystyle{abbrv}
\bibliography{bibfile-4}

\end{document}